\documentclass[11pt, reqno]{amsart}
\usepackage[utf8]{inputenc}
\usepackage{amsmath}
\usepackage{amsthm}
\usepackage{hyperref}
\usepackage{amssymb}
\usepackage{amsfonts}
\usepackage{caption}
\usepackage{MnSymbol}
\usepackage[new]{old-arrows}
\usepackage{circuitikz, booktabs}
\usepackage{tcolorbox}
\usepackage{textcomp}
\usepackage{orcidlink}
\usepackage{relsize}
\usepackage{MnSymbol}
\usepackage{extpfeil}
\usepackage{blindtext,amsmath,amsthm,amsfonts, textcomp, mathrsfs, mathtools}
\usepackage[shortlabels]{enumitem}

\usepackage[top=25mm,bottom=25mm,left=27mm,right=27mm]{geometry}
\newtheorem{theorem}{Theorem}[section]

\newtheorem{corollary}[theorem]{Corollary}

\newtheorem*{theorem*}{Theorem}

\newtheorem*{remark*}{Remark}
\newtheorem*{problem*}{Problem}
\newtheorem*{conjecture*}{Conjecture}
\newtheorem*{question*}{Question}
\newtheorem{lemma}[theorem]{Lemma}
\newtheorem{proposition}[theorem]{Proposition}

\newcommand{\rom}[1]{\uppercase\expandafter{\romannumeral #1\relax}}
\newcommand{\Q}{\mathbb{Q}}
\newcommand{\Z}{\mathbb{Z}}
\newcommand{\N}{\mathbb{N}}
\newcommand{\rL}{\mathcal{L}}

\newcommand{\rK}{\mathcal{K}}
\newcommand{\rN}{\mathcal{N}}

\newcommand{\R}{\mathbb{R}}
\newcommand{\C}{\mathbb{C}}
\newcommand{\F}{\mathbb{F}}

\def\house#1{{%
    \setbox0=\hbox{$#1$}
    \vrule height \dimexpr\ht0+1.4pt width .4pt depth \dp0\relax
    \vrule height \dimexpr\ht0+1.4pt width \dimexpr\wd0+2pt depth \dimexpr-\ht0-1pt\relax
    \llap{$#1$\kern1pt}
    \vrule height \dimexpr\ht0+1.4pt width .4pt depth \dp0\relax
}}

\title[Bogomolov Property for Abelian Varieties]{Lower bound for the canonical height on Abelian varieties over totally $p$-adic extensions}

\author[S.Kala]{Sushant Kala\, \orcidlink{0009-0004-7689-0131}}
\address[]{Chennai Mathematical Institute, H1, SIPCOT IT Park, Kelambakkam, Siruseri, Tamil Nadu 603103, India}
\email{sushantkala786@gmail.com, sushantkala@cmi.ac.in}


\begin{document}

\maketitle
 \textbf{Abstract.} Let $A/\mathbb{Q}$ be an abelian variety and let $\hat{h}$ be the canonical height on $A(\overline{\Q})$ associated to a symmetric ample line bundle $\mathcal{L}$ on $A$. We prove that $\hat h$ is bounded away from zero on non-torsion points of $A$ defined over the maximal totally $p$-adic extension of $\mathbb{Q}$, for all but finitely many primes $p$. More generally, for abelian varieties over a number field $K$, we obtain a similar height gap over certain infinite extensions of $K$, including Galois extensions with finite local degree at non-archimedean places.

 % In this article, we prove that the N\'{e}ron–Tate height of non-torsion points in $A(\mathcal{K})$ is bounded below by an absolute constant depending only on $A$, $K$ and $\mathcal{K}$. As a consequence, we obtain the Bogomolov property for totally $p$-adic points of an abelian variety $A/\Q$ in the case when $A$ has good reduction at $p$.  Moreover, our result also implies the finiteness of torsion points in $A(\mathcal{K})$.

 \bigskip

 % In this paper, we provide an explicit lower bound for the canonical height of non-torsion points of abelian varieties lying in an ASP. As a consequence, we obtain the effective Bogomolov property for totally $p$-adic points of an abelian variety defined over the rationals, where $A$ has good reduction at $p$.

\section{\bf Introduction}

\bigskip

The study of points of small height is a classical theme in Diophantine geometry. In 2001, E. Bombieri and U. Zannier \cite{BZ} introduced the notions of the Bogomolov (B) and Northcott (N) properties for subsets of algebraic numbers with respect to the logarithmic Weil height. Let $h$ denote the logarithmic Weil height and let $S \subseteq \overline{\Q}$. We say that $S$ satisfies the Bogomolov property (B) if
$$
\liminf_{\alpha \in S \setminus \mu_\infty} h(\alpha) > 0,
$$
where $\mu_\infty$ denotes the set of all roots of unity. More precisely, this property characterizes subsets of algebraic numbers that do not contain points of arbitrarily small positive height. The terminology ``Bogomolov property" is inspired by the celebrated Bogomolov conjecture, which concerns a similar phenomenon in a more general geometric setting. Further, we say $S$ satisfies the Northcott property (N) if 
$$
\liminf_{\alpha \in S} h(\alpha) = \infty.
$$
It follows from a theorem of Northcott that subsets of algebraic numbers of bounded degree satisfy the Northcott property, and hence also the Bogomolov property. It is therefore natural to ask when these properties hold for infinite extensions of $\Q$. The same question has also been studied in a geometric setting for elliptic curves and abelian varieties with respect to canonical height.\\

Let \(A/K\) be an abelian variety defined over a number field \(K\). Let \(\mathcal{L}\) be a symmetric ample line bundle on \(A\) and \(\hat{h}\) be the associated canonical height on $A(\overline{K})$. For an infinite extension $\mathcal{K}/K$, we say that $A(\mathcal{K})$ satisfies the Bogomolov property if

$$
\liminf_{P \in A(\mathcal{K})\setminus A_{\mathrm{tor}}(\mathcal{K})} \hat{h}(P) > 0,
$$
where $A(\mathcal{K})$ denotes the set of points of $A(\overline{K})$ defined over $\mathcal{K}$, and $A_{\mathrm{tor}}(\mathcal{K})$ denotes the set of torsion points in $A(\mathcal{K})$. Therefore, a natural question is to determine for which infinite extensions $\mathcal{K}/K$ the set $A(\mathcal{K})$ satisfies the Bogomolov property.\\

The first such result was obtained by S. Zhang \cite{Zhang} in 1998. Using an equidistribution theorem of Szpiro--Ullmo--Zhang \cite{SUZ}, Zhang showed that if $A/\Q$ is an abelian variety and $\Q^{tr}$ denotes the maximal totally real extension of $\Q$, then $A(\Q^{tr})$ satisfies the Bogomolov property. In 2004, Baker and Silverman \cite{Baker-Silverman} proved that, for any abelian variety $A/K$, the set $A(K^{ab})$ satisfies the Bogomolov property, where $K^{ab}$ denotes the maximal abelian extension of $K$ in $\overline{K}$.\\

The non-archimedean analogue of Zhang's theorem for elliptic curves was established by M.~Baker and C.~Petsche \cite{Baker-2}. Let $v$ be a finite place of $K$. An infinite Galois extension $K'/K$ is said to be totally $v$-adic of type $(e,f)$ if there exists an embedding $K' \hookrightarrow K_v$, where $e$ and $f$ are the corresponding ramification index and inertia degree. Baker and Petsche showed that if $E/K$ is an elliptic curve and $\mathcal{K}/K$ is a totally $v$-adic extension of type $(e,f)$, then $E(\mathcal{K})$ satisfies the Bogomolov property. This was later extended to asymptotically positive extensions by Dixit and the author in \cite{AB-SK}.\\

For higher-dimensional abelian varieties, the first non-archimedean analogue of Zhang's theorem was obtained by W.~Gubler \cite{Gubler} using methods from tropical geometry. Let $K^{nr,v}$ denote the maximal algebraic extension of $K$ unramified at $v$. If $A/K$ is totally degenerate at $v$, then Gubler showed that $A(K^{nr,v})$ satisfies property (B). Since every totally $v$-adic field of type $(e,f)$ can be realized inside $K^{nr,v}$ after a suitable finite extension of the base field, this implies the Bogomolov property for totally $v$-adic points in the totally degenerate case. It follows from the recent work of A.~Plessis and S.~Sahoo \cite{Plessis} that the totally degenerate hypothesis can be relaxed to bad reduction.\\

However, much less is known when $A/K$ has good reduction at a place $v$ of $K$, even in the case $K=\Q$. For an abelian variety $A/\Q$ with good reduction at $p$, the best known result is due to A. Chambert-Loir \cite{CL}, who showed that the Bogomolov property holds for a Zariski dense open subset of $A(\Q^{tp})$. In this article, we establish property (B) for $A(\mathcal{K})$, where $\mathcal{K}/K$ is an infinite extension satisfying certain asymptotic splitting conditions at a non-archimedean place of $K$. In particular, our result implies the Bogomolov property for $A(\Q^{tp})$ for all but finitely many places of good reduction. To state our result, we recall the notion of asymptotically positive extension, introduced by A. Dixit and the author in \cite{AB-SK}. \\

For a number field $F$ and a rational prime power $q=p^f$, let $\rN_q(F)$ denote the number of prime ideals of $F$ of norm $q$. An infinite algebraic extension $\rK/\Q$ can be written as a tower
$$
\rK \supsetneq \cdots \supsetneq K_m \supsetneq K_{m-1} \supsetneq \cdots \supsetneq K_0=\Q,
$$
where each $K_i$ is a number field. Define
$$
\psi_q(\rK) := \lim_{i\to \infty} \frac{\rN_q(K_i)}{[K_i:\Q]}.
$$
This limit exists and is well-defined, i.e. it is independent of the choice of the tower ${K_i}$ (see \cite[Proposition 4.2]{AB-SK}). When the extension $\rK$ is clear from the context, we simply write $\psi_q$ for $\psi_q(\rK)$. We say that $\rK$ is an \textit{asymptotically positive extension} if $\psi_q > 0$ for some prime power $q$. For example, the maximal totally $p$-adic extension $\Q^{tp}$ of $\Q$ is an asymptotically positive extension with $\psi_p = 1$. \\

The main result of this article is as follows.

\begin{theorem}\label{Main 2}
     Let $A/K$ be a geometrically simple abelian variety of dimension $g$, and let $\hat{h}$ be the canonical height associated with a very ample line bundle $\mathcal{L}$ on $A$. Then there exist a finite set $\mathcal{P}$ of rational primes and an absolute constant $C(A,K, \mathcal{\rL})>0$, depending only on $A$, $\rL$ and $K$, such that for every infinite algebraic extension $\mathcal{K}/K$,
$$
\liminf_{P\in A(\mathcal{K})} \hat{h}(P)
\geq
C(A, K , \rL) \sum_{q=p^k}' \psi_q(\mathcal{K})\frac{\log q}{q^{2g}},
$$
where the sum runs over all prime powers $q=p^k$ with $p\notin \mathcal{P}$. 
\end{theorem}

\medskip
\noindent
\textbf{Remark 1.} As a consequence of Theorem \ref{Main 2}, we obtain the Bogomolov property for $A(\mathcal{K})$ whenever $\psi_{p^k}(\mathcal{K})>0$ for some prime $p\notin\mathcal{P}$. It is important to note that the lower bound in Theorem \ref{Main 2} combines the contributions from different primes linearly. In the classical setting of the Weil height, lower bounds of this nature were first obtained by Bombieri and Zannier \cite{BZ}, and later by P. Fili \cite{Fili1}. For elliptic curves, such lower bounds were established by Baker and Petsche \cite{Baker-2}. In \cite{AB-SK}, Dixit and the author established analogous lower bounds both in the classical setting and for elliptic curves. To the best of our knowledge, Theorem \ref{Main 2} provides the first such result in higher dimensions.\\

As an immediate consequence of Theorem \ref{Main 2}, we obtain the following corollary.

\medskip

\begin{corollary}\label{Bogomolov_totally_p-adic}
  In the setting of Theorem 1, suppose that $\mathcal{K}/K$ is an infinite Galois extension having finite local degree at a non-archimedean place $v$ of $K$ whose residue characteristic does not belong to $\mathcal{P}$. Then
$$
\liminf_{P\in A(\mathcal{K})}
\hat h(P)
\geq
C(A,K,\mathcal L)\frac{\log q_v}{q_v^{2g}}.
$$
\end{corollary}

\medskip

As a special case, take $K=\Q$ and $\rK = \Q^{(tp)}$, the maximal algebraic extension of $\Q$ in which $p$ splits completely.

\medskip

\begin{corollary}
   Let $A/\Q$ be a geometrically simple abelian variety of dimension $g$, and let $\hat{h}$ be the canonical height associated with a very ample line bundle $\mathcal{L}$ on $A$.  Then there exist a finite set $\mathcal{P}$ of rational primes and an absolute constant $C(A, \rL) > 0$, depending only on $A$ and $\rL$, such that for every $p \notin \mathcal{P}$,
    $$
    \liminf_{P \in A(\Q^{tp})} \hat{h}(P) \geq C(A, \rL) \frac{\log p}{p^{2g}}.
    $$
\end{corollary}

\medskip

\noindent
\textbf{Remark 2.} The corollary above establishes the Bogomolov property for $A(\Q^{tp})$ for all but finitely many primes $p$ when $A/\Q$ is geometrically simple. In particular, it applies to all but finitely many primes of good reduction for $A$. In higher dimensions ($g \geq 2$), the Bogomolov property for $A(\Q^{tp})$ was previously known only for primes at which $A$ has bad reduction (see \cite{Gubler, Plessis}). Thus, even in the geometrically simple case, Corollary \ref{Bogomolov_totally_p-adic} provides the first examples of primes of good reduction $p$ for which $A(\Q^{tp})$ satisfies the Bogomolov property.\\

\noindent
\textbf{Remark 3.} In all the results mentioned above, the liminf is taken over all points in $A(\mathcal{K})$. Consequently, whenever $\psi_{p^k}(\mathcal{K})>0$ for some $p\notin\mathcal{P}$, we obtain both the Bogomolov property for $A(\mathcal{K})$ and the finiteness of the torsion subgroup $A_{\mathrm{tor}}(\mathcal{K})$.

\medskip

The hypotheses in Theorem \ref{Main 2} that $A/K$ is geometrically simple and that $\mathcal{L}$ is very ample can be relaxed, at the cost of a slightly weaker lower bound. In the general setting, we obtain the following result.

\begin{theorem}\label{Main 3}
      Let $A/K$ be an abelian variety of dimension $g$, and let $\hat{h}$ be the canonical height associated with a symmetric ample line bundle $\mathcal{L}$ on $A$. Then there exists a finite set $\mathcal{P}$ of rational primes, a finite extension $K'/K$ and an absolute constant $C(A,K', \mathcal{\rL})>0$, depending only on $A$, $\rL$ and $K'$, such that for every infinite algebraic extension $\mathcal{K}/K$,
$$
\liminf_{P\in A(\mathcal{K}K')} \hat{h}(P)
\geq
C(A,K', \mathcal{L})\sum_{q=p^k}' \psi_q(\mathcal{K}K')\frac{\log q}{q^{2g}},
$$
where the sum runs over all prime powers $q=p^k$ with $p\notin \mathcal{P}$. 
\end{theorem}

\medskip
\noindent
\textbf{Remark 4.} Suppose that $\mathcal{K}/K$ is an asymptotically positive extension with $\psi_{p^k}(\mathcal{K})>0$ for some prime power $p^k$, where $p\notin\mathcal{P}$. A priori, the Bogomolov property for $A(\mathcal{K})$ does not follow immediately from Theorem \ref{Main 3}, since it may happen that $\psi_{p^k}(\mathcal{K}K')=0$. However, since asymptotic positivity is preserved under finite extensions, we have $\psi_{p^{k'}}(\mathcal{K}K')>0$ for the same prime $p$ and possibly some $k'\neq k$ (see \cite[Lemma 4.4]{AB-SK}). Hence $A(\mathcal{K}K')$ satisfies the Bogomolov property, and therefore so does $A(\mathcal{K})$.\\

\noindent
\textbf{Remark 5.} In the special case, when $K=\Q$ and $\mathcal{K}=\Q^{(tp)}$, Theorem \ref{Main 3} implies that if $A/\Q$ is an abelian variety and $\hat{h}$ is the canonical height associated with an ample line bundle on $A$, then $A(\Q^{(tp)})$ satisfies the Bogomolov property for all but finitely many primes $p$. Thus, for all but finitely many primes $p$, this gives an affirmative answer to a question of V. Dimitrov \cite{Dimitrov}.\\

The main idea in the proof of Theorem \ref{Main 2} is the construction of an auxiliary function together with zero-estimate techniques from transcendental number theory in a geometric setting, as used by S. David and M. Hindry \cite{David-Hindry} in the archimedean case. The new ingredient is to adapt this method to the non-archimedean setting.\\

The organization of this article is as follows.

% The key ideas in the proof of Theorem \ref{Main 1} and Theorem \ref{Main 2} are inspired by the classical auxiliary function and zero estimate techniques from transcendental number theory in the geometric setting, as employed by David--Hindry \cite{David-Hindry}. 

\subsection{Organization} There are 7 sections including the introduction. In Section 2, we recall key definitions and results concerning canonical heights and discuss the notion of differentiating sections on an abelian variety. In Section 3, we construct an auxiliary section $F$ vanishing at the origin using Siegel's lemma. Section 4 is devoted to height estimates for this auxiliary section and the choice of parameters. The set of rational primes $\mathcal{P}$ appearing in the statements of our theorems is also chosen in Section 4 and remains fixed throughout the paper. In Section 5, we apply Philippon's zero estimate to bound the number of small-height points at which $F$ vanishes.
Finally, Sections 6 and 7 prove Theorems \ref{Main 2} and \ref{Main 3}, respectively.

\section{\bf Preliminaries}

\subsection{Heights} A place $v$ of a number field $K$ is an absolute value $\left| \ \cdot \ \right|_v: K \rightarrow [0, \infty)$ whose restriction to $\mathbb{Q}$ is either the standard absolute value or a $p$-adic absolute value for some prime $p$ with $|p|_v=p^{-1}$. We set $d_v=\left[K_v: \mathbb{R}\right]$ in the former and $d_v=\left[K_v: \mathbb{Q}_p\right]$ in the latter case. The absolute Weil height of a point $\textbf{x}=\left[x_0: \ldots: x_n\right] \in \mathbb{P}_{\mathbb{Q}}^n(K)$ with $x_0, \ldots, x_n \in K$ is

$$
h(\textbf{x})=\frac{1}{[K: \mathbb{Q}]} \sum_{v \in M_K} d_v \log \max \left\{\left|x_0\right|_v, \ldots,\left|x_n\right|_v\right\},
$$

\noindent
where $M_K$ denotes the set of all places $v$ of $K$. The factor $d_v$ on the right-hand side ensures that the normalized absolute values $\left| \ \cdot \ \right|_v^{d_v}$ satisfy the product formula. This implies $h(\mathbf{x})$ is independent of the choice of projective coordinates. Moreover, if we consider another number field containing the coordinates of $\mathbf{x}$, the Weil height remains unchanged. Therefore, $h(\cdot)$ is well-defined on $\mathbb{P}^n_{\mathbb{Q}}(\bar{\Q})$, where $\bar{\Q}$ is a fixed algebraic closure of $\Q$. For a comprehensive account of the Weil height and its properties, we refer to \cite{ Bombieri-Gubler, Hindry-Silverman}.\\

\subsection{Arakelov height} Let $\textbf{x} = [x_0 : x_1 : \ldots : x_n] \in \mathbb{P}^n(\overline{\Q})$ such that $\textbf{x}$ is defined over a number field $K$. The Arakelov height of $\textbf{x}$ is defined as
$$
h_{Ar}(\textbf{x})=\frac{1}{[K: \Q]}\left(\sum_{v \in M_{K} \setminus M_{K}^\infty} d_{v} \log \max _{0 \leq i \leq n}\left|x_{i}\right|_{v}+\sum_{v \in M_{K}^{\infty}} d_{v} \log \sqrt{\sum_{0 \leq i \leq n}\left|x_{i}\right|_{v}^{2}}\right),
$$
where $M_{K}^{\infty}$ denote the set of archimedean places of $K$. This new height differ from the Weil height only in using the $L^2$ local height instead of $L^{\infty}$, at the infinite places. Let $W$ be an $m$-dimensional subspace of $\overline{\Q}^n$. The $m$-th exterior power $\bigwedge^m W$ is a one-dimensional subspace of $\bigwedge^m \overline{\Q}^n$. Therefore, we may view $W$ as a point $P_W$ of the projective space
$\mathbb{P}(\bigwedge^m \overline{\Q}^n)$. The latter may be identiﬁed with projective space of dimension $\binom{n}{m}$ by using standard coordinates.\\

\noindent
The Arakelov height of $W$ is defined as 
$h_{Ar}(W) : = h_{Ar}(P_W)$. Moreover, for a matrix $A$ with coefficients in $\overline{\Q}$, we define 
$$
h_{Ar}(A): = h_{Ar}(\text{row}(A)),
$$
where $\text{row}(A)$ is the row space of $A$ as a $\overline{\Q}$-vector space.

% \begin{lemma}\label{Orthogonal}
%     If $W$ be $\overline{\Q}$-vector subspace of $\overline{\Q}^n$, then
%     $$
%     h_{Ar}(W)= h_{Ar}(W^\perp),
%     $$
%     where $W^\perp$ is the orthogonal complement of $W$ in $\overline{\Q}^n$.
% \end{lemma}
% \begin{proof}
%     See \cite[Chapter I, Lemma 4C]{Schmidt}.
% \end{proof}

\subsection{N\'{e}ron-Tate height of abelian varities} Let $A/K$ be an abelian variety defined over a number field $K$ and $\rL$ a symmetric line bundle on $A/K$. We fix an embedding $\psi : A \rightarrow \mathbb{P}_K^n$ with $\rL = \psi^{*}(\mathcal{O}(1))$. Let $h_{A, \, \mathcal{L}}$ be the height function on $A(\overline{K})$ corresponding to the Weil height on $\mathbb{P}_K^n$ under the embedding $\psi$. It is easy to show that $h_{A, \mathcal{L}}$ satisfies the finiteness property (i.e., Northcott's theorem) and hence provides a partial order on $A(L)$ for finite extensions $L/K$. Despite this, $h_{A, \mathcal{L}}$ is not compatible with the group structure of an abelian variety. N\'{e}ron and Tate introduced the canonical height on $A(\overline{K})$, which is a quadratic form on $A(\overline{K}) \otimes \R$ and is unique up to the linear equivalence of line bundles.

\medskip

\begin{theorem}{(N\'{e}ron-Tate)} Let $A/K$ be an abelian variety and $\mathcal{L}$ a symmetric ample line bundle on $A/K$. Then there exists a unique function, called the canonical height on $A(\overline{K})$ relative to $\mathcal{L}$, 
$$
\widehat{h}_{A, \, \mathcal{L}} \, : \, A(\overline{K}) \longrightarrow \mathbb{R},  
$$
such that 
\begin{enumerate}[(i)]
    \item $\widehat{h}_{A, \, \mathcal{L}}(P) = h_{A, \, \mathcal{L}}(P) + O(1)$ for all $P \in A(\overline{K})$,

    \item $\widehat{h}_{A, \, \mathcal{L}}([2] \cdot P) = 4 \, \widehat{h}_{A, \, \mathcal{L}}(P)$ for all $P \in A(\overline{K})$.
\end{enumerate}
\end{theorem}

\medskip

% The canonical height function satisfies various interesting properties and can be expressed as a limit of Weil height as we see in the following proposition.

The N\'{e}ron-Tate height depends on the choice of the symmetric ample line bundle. We often use $\hat{h}$ to denote the canonical height on abelian variety coming from a symmetric ample line bundle $\rL$ when there is no ambiguity. Tate's limit process induces the N\'{e}ron-Tate or canonical height on $A(\bar{K})$. If $K$ is a number field, the N\'{e}ron-Tate height of a point $P \in A(\overline{K})$ is defined as
$$
\hat{h}(P)=\lim _{N \rightarrow \infty} \frac{h\left(\left[2^N\right]\cdot P\right)}{4^N}.
$$
\noindent
The N\'{e}ron-Tate height function on an abelian variety $A/K$ has many interesting properties.

\medskip

\begin{proposition} \label{CHAB}
    Let $A/K$ be an abelian variety and $\hat{h}$ the canonical height on $A$ with respect to a symmetric ample line bundle $\mathcal{L}$. Then, for $P, Q \in A(\overline{K})$
    \begin{enumerate}[(i)]
        % \item 
        % $\hat{h}(P) = \lim_{n \rightarrow \infty} \frac{1}{4^n} h_{A, \, \mathcal{L}}([2^n] (P)).$
        \item 
        $\hat{h}(P) = 0$ iff $P \in A_{tor}$,
        \item
        $\hat{h}(P+Q) + \hat{h}(P-Q) = 2\, \hat{h}(P) + 2 \, \hat{h}(Q)$,
        \item 
        $\hat{h}([n] \, P) = n^2 \, \hat{h}(P)$ for all $n \in \Z$.
    \end{enumerate}
\end{proposition}

\noindent
Interested readers can find proofs of the above properties and more detailed discussions on canonical heights on abelian varieties in \cite[Section B.5]{Hindry-Silverman}.
\medskip

% Among them, we know that $\hat{h}_A$ defines a positive definite quadratic
% form on $A(\overline{K}) \otimes_{\Z} \R$. Moreover, it classifies the torsions points of $A(\overline{K})$ precisely as the points of height zero.\\

\subsection{Asymptotically positive extensions}
We recall here the definition of asymptotically positive extensions introduced by A. Dixit and the author in \cite{AB-SK}. For a number field $K$ and a rational prime power $q=p^f$, define
\begin{equation*}
    \rN_q(K) := \text{ the number of prime ideals of } K \text{ with norm } q.
\end{equation*}
An infinite extension $\rK/\Q$ can be written as a tower of number fields
$$
\rK \supsetneq \cdots \supsetneq K_m \supsetneq K_{m-1} \supsetneq \cdots \supsetneq K_0 = K,
$$
where $K_i/\Q$ are finite extensions. Define
$$
    \psi_q(\rK) = \psi_q := \lim_{i\to \infty} \frac{\mathcal{N}_q(K_i)}{[K_i:\Q]}.
$$
This limit exists and is well-defined, i.e., independent of the tower $\{K_i\}$ (see \cite[Proposition 4.2]{AB-SK}) and $0\leq \psi_q \leq 1$. For instance, if a prime $p$ splits completely in $\rK$, then $\psi_p =1$ and $\psi_{p^f} =0$ for all $f>1$. If $\psi_q > 0$ for some prime power $q$, we call the extension $\rK$ to be \textit{asymptotically positive}. We refer to \cite[Section 4]{AB-SK} for a comprehensive account of these extensions.\\

For a finite extension $L/K$ of number fields, let $S$ be a finite subset of $M_L$. We need the following well-known inequality from Diophantine approximation.

\begin{lemma}(Liouville's inequality)
  If $\alpha \in L$ and $\beta \in K$ with $\alpha \neq \beta$, then
\begin{equation}\label{Liouville's inequality}
    (2 H(\alpha) H(\beta))^{-[L: K]} \leq \prod_{w \in S}\|\alpha-\beta\|_{w, K} \leq(2 H(\alpha) H(\beta))^{[L: K]},
\end{equation}
where $\|x\|_{w, K}=\left|N_{L_w / K_v}(x)\right|_v$.
\end{lemma}
\begin{proof}
    See \cite[Theorem 1.5.21]{Bombieri-Gubler}.
\end{proof}

\subsection{Preliminaries on differential operators}

% We briefly review here the notations and results of [16], section 4.1. 
Let $A/K$ be an abelian variety of dimension $g$. Under the embedding $\psi : A \longrightarrow \mathbb{P}_K^n$ associated to a very ample line bundle $\mathcal{L}$, let $\mathfrak{U}$ be the prime ideal of definition of $A$ in $K\left[X_0, \ldots, X_n\right]$, and let $\mathcal{A}$ be the coordinate ring $K\left[X_0, \ldots, X_n\right] / \mathfrak{U}$. Let $u \in A(K)$ and $ \mathfrak{m}$ be its ideal of definition in $A$. Denote by $ \mathcal{A}_{\mathfrak{m}}$ the local ring at the point $u$ and $\hat{\mathcal{A}}_{\mathfrak{m}}$ its completion with respect to the $\mathfrak{m}$-adic topology. It is well known that $\hat{\mathcal{A}}_{\mathfrak{m}}$ is isomorphic to $K\left[\left[T_1, \ldots, T_g\right]\right]$. In particular, when $u = \mathcal{O}$ is the identity element of $A(K)$, we fix such an isomorphism and denote it by $\Phi_0$. If we denote by $s : A \times A \rightarrow A$ the addition in $A$ and $s^*$ its inverse image at the level of $K$-algebras, we can define for all $u \in A(K)$ isomorphisms $\Phi_u$ of $\hat{A}_m$ with $K\left[\left[T_1, \ldots, T_g\right]\right]$ by

$$
\Phi_u=\Phi_0 \circ(\mathrm{Id} \otimes \pi) \circ s^*: \hat{\mathcal{A}}_{\mathfrak{m}} \xrightarrow{\sim} \hat{\mathcal{A}}_{\mathfrak{m}_0} \otimes_K\left(\hat{\mathcal{A}}_{\mathfrak{m}} / \mathfrak{m} \hat{\mathcal{A}}_{\mathfrak{m}}\right) \simeq \hat{\mathcal{A}}_{\mathfrak{m}_0}
$$
\noindent
where $\pi$ denotes the canonical projection $\pi: \hat{\mathcal{A}}_{\mathfrak{m}} \rightarrow \hat{\mathcal{A}}_{\mathfrak{m}} / \mathfrak{m} \hat{\mathcal{A}}_{\mathfrak{m}} \simeq K$ and $m_0$ is the ideal corresponding to $\mathcal{O}$.
Now, we define a basis for the space $\operatorname{Hom}_K\left(\hat{\mathcal{A}}_{\mathfrak{m}}, K\right)$ using the family of differential operators $\left(\partial_{\Phi_u}^k\right)_{k \in \mathbb{N}^g}$, by setting
$$
\partial_{\Phi_u}^k=\left.\frac{1}{k!} \frac{\partial^{|k|}}{\partial \mathbf{T}^k}\right|_{\mathbf{T}=0},
$$
with $k=\left(k_1, \ldots, k_g\right)$, $|k|=k_1+\ldots+k_g$, and $k! = k_{1}! \ldots k_{g}!$ and $\mathbf{T}^k = T_1^{k_1} T_2^{k_2} \ldots T_g^{k_g} $. Let $\epsilon_1, \ldots, \epsilon_g$ be the vectors of the canonical basis of $\R^g$, then the $\partial_{\Phi_u}^i$ form a basis of the tangent space $\operatorname{Hom}\left(\mathfrak{m} / \mathfrak{m}^2, K\right)$ of $A$ at $u$.

\medskip

We now define what it means for a homogeneous polynomial in $K[X_0,\ldots,X_n]$ to vanish to order at least $T+1$ at a point $u \in A$. Let $P$ be such a polynomial of degree $D$, let $u$ be a point of $A$, and let $F$ be a homogeneous form in $K[X_0,\ldots,X_n]$ of degree $D$ that does not vanish at $u$. Denote by $P_{u,F}$ the image of $\frac{P}{F}$ in $\hat{\mathcal{A}}_{\mathfrak{m}}$. We say that $P$ vanishes to order $\geq T+1$ at $u$ if, for every affine neighborhood of $u$ of the form $\{F \neq 0\}$, we have
$$
\partial_{\Phi_u}^k P_{u,F} = 0 \quad \text{for all } k \in \mathbb{N}^g \text{ with } |k| \leq T.
$$
\medskip
% We also define for $\mathcal{V}$ a vector subspace of $T_{A(\C)}$ of dimension $d$, the notion of vanishing $P$ at the point $u$ to order $\geq T+1$ along $\mathcal{V}$. To do this, we choose a basis of $\mathcal{V}$ that we complete into a basis of $T_{A(\C)}$ using the $\partial_{\Phi_0}^{\epsilon_i}$ and we agree that these vectors are numbered from $d+1$ to $g$. We say that $P$ is zero at $u$ to order $\geq T+1$ along $\mathcal{V}$ if $\partial_{\Phi_u}^k P_{u, F}=0$ for all $k \in \N^d$ such that $|k| \leq T$ and where we identify $k \in \N^d$ with the element of $\N^g$ whose last $g-d$ coordinates are zero. 
Let $\alpha$ be an admissible isogeny on $A$ of degree $q(\alpha)$, denote by $\alpha^*$ the algebra morphism deduced from $\alpha$ given by
\begin{align*}
\alpha^*: \mathcal{A} & \longrightarrow \mathcal{A} \\
P & \longmapsto P\left(P_{0, \alpha}, \ldots, P_{n, \alpha}\right),
\end{align*}
where ($P_{0, \alpha}, \ldots, P_{n, \alpha}$ ) is a family of homogeneous polynomials of degree $q(\alpha)$ representing the action of $\alpha$ on $A$. Let $u \in A$ and $s_u^*: \mathcal{A} \rightarrow \mathcal{A} \otimes \hat{\mathcal{A}}_m$ be the algebra morphism deduced from addition in $A$ in the neighborhood of $(0, u)$.

\medskip
\noindent
\textbf{Definition.} Let $u \in A, \ \alpha$ be an admissible isogeny and $k \in \N^g$. For all $P \in \mathcal{A}$, we note $\partial_{u, \alpha}^k(P)$ the coefficient of the monomial $\textbf{T}^k$ of
$$
\left(\mathrm{Id} \otimes \Phi_u\right) \circ \alpha^* \otimes \mathrm{Id} \circ s_u^*(P) \in \mathcal{A} \otimes K\left[\left[T_1, \ldots, T_g\right]\right].
$$

The Taylor series $\sum_{k \in \N^n} \partial_{u, \alpha}^k(P) \mathbf{T}^k$ is the image of $P$ under the map:
$$
\mathcal{A} \xrightarrow{s_{\mathrm{u}}^*} \mathcal{A} \otimes \hat{\mathcal{A}}_{\mathfrak{m}} \xrightarrow{\alpha^* \otimes \mathrm{Id}} \mathcal{A} \otimes \hat{\mathcal{A}}_{\mathfrak{m}} \xrightarrow{\mathrm{Id} \otimes \Phi_{\mathrm{u}}} \mathcal{A} \otimes K\left[\left[T_1, \ldots, T_g\right]\right] .
$$

\noindent
It is easy to verify that $\partial_{u, \alpha}^k(P) \in \mathcal{A}$ and, in the vicinity of the origin, the nullity of $\partial_{u, \alpha}^k(P)$ does not depend on the choice of forms representing addition or isogeny $\alpha$.

% Let $k$ be a positive integer, $\mathcal{V}$ be a vector subspace of $T_{A(\C)}$ of dimension $d$, and $\mathfrak{I}$ be a homogeneous ideal of $\mathcal{A}$. We denote by $\partial_{u, \alpha}^{k, \mathcal{V}} \mathfrak{J}$ the ideal

% $$
% \left(\partial_{u, \alpha}^{k, \mathcal{V}}(Q), \ k \in \N^d, \ |k| \leq k, \ Q \in \mathfrak{J}\right),
% $$

% \noindent
% where we consider $k \in \N^g$ and the first $d$ directions of derivation correspond to a basis of $\mathcal{V}$.

% The following lemma connects the notations introduced above:

% \begin{lemma}
% Let $u, h \in A, \ \alpha$ be an admissible isogeny, $k \in \N, \mathfrak{J}$ be a homogeneous ideal of $\mathcal{A}$, and $\mathcal{V}$ be a $d$-dimensional vector subspace of $T_{A(\C)}$. The following properties are equivalent:
% \begin{enumerate}[(i).]
%     \item $h \in \mathcal{Z}\left(\partial_{u, \alpha}^{k, \mathcal{V}} \mathfrak{J}\right)$.
%     \item Any polynomial of $\mathfrak{J}$ vanishes to order $\geq k+1$ at $u+\alpha(h)$ along $\mathcal{V}$.
% \end{enumerate}
% \end{lemma}

\subsection{Embedding $A$ into $A \times A$} We denote by $\mathcal{M}$ the line bundle $\mathcal{L} \otimes \mathcal{L}$ on $A \times A$ associated to $\mathcal{L}$. For an integer $N$, denote by $N$ the endomorphism of $A$ given by multiplication by $N$. We have
\begin{align*}
 A & \stackrel{i}{\longrightarrow} A \times A \ \stackrel{\rL \otimes \rL} \longrightarrow \mathbb{P}^{n} \times \mathbb{P}^{n} \stackrel{\text { Segre }}{\longrightarrow} \mathbb{P}^{(n+1)^{2}-1}
\end{align*}
given by $P \stackrel{i}\longmapsto(P,[N] \cdot P)$.

\medskip

 Denote by $\Gamma^{0}(A \times A, \mathcal{M})$ the space of all global sections of the line bundle $\mathcal{M}$ on the variety $A \times A$. Let $\left\{s_{0}, \ldots, s_{l}\right\}$ be a basis of $\Gamma^{0}(A \times A, \mathcal{M})$. We fix an integer $L$ to be chosen later. By projective normality, one can choose a basis $\left\{f_{1}, \ldots, f_{m}\right\}$ of the $K$-vector space $\Gamma^{0}\left(A \times A, \mathcal{M}^{\otimes L}\right)$ such that all $f_{i}$'s are homogeneous of degree $L$ in the $s_{i}$'s. Furthermore, the $s_{i}$'s may also be seen as (1,1)-homogeneous forms of $K[\mathbf{X}, \mathbf{Y}]$ where $\mathbf{X}=\left(X_{0}, \ldots, X_{n}\right)$ and $\mathbf{Y}=\left(Y_{0}, \ldots, Y_{n}\right)$. Finally, we denote by $T_A$ the tangent space at the origin of the abelian subvariety $B=i(A)$ of $A \times A$ defined by $y=[N] x$.

\medskip

% Consider the embedding

% $$
% \begin{aligned}
% \iota: A & \longrightarrow A \times A \\
% P & \longmapsto(P, [N] \cdot P)
% \end{aligned}
% $$

% \noindent
% and $\varphi_{\mathcal{L}}$ a projectively normal embedding associated with the bundle $\mathcal{L}=\mathcal{L}_0^{\otimes 4}$. Let $\mathcal{M}=$ $\pi_1^* \mathcal{L} \otimes \pi_2^* \mathcal{L}$ be the bundle on $A \times A$ constructed using the projections onto the first and second factors. 

\section{\bf Siegel lemma and section with higher order vanishing at $0$}

Our goal is to construct an auxiliary section $F$ of $\Gamma^0\left(A \times A, \mathcal{M}^{\otimes L}\right)$ with higher order of vanishing at the origin.
% Let $s_0, \ldots, Z$ be a basis of the $K$-vector space $\Gamma\left(A \times A, \mathcal{M}\right)$ and $f_1, \ldots, f_m$ be monomials of degree $L$ in the $s_i$, forming a basis of $\Gamma\left(A \times A, \mathcal{M}^{\otimes L}\right)$. 
Denote by $L_w$ the dimension of the $K$-vector space $\Gamma^0\left(A \times A, \mathcal{M}^{\otimes L}\right)$. The section $F$ will therefore be of the form
$$
F=\sum_{i=1}^m b_i f_i,
$$
\noindent
where the coefficients $b_i$ are in $K$. Using the earlier notations, we fix a basis $\partial_{\Phi_0}^{\epsilon_1}, \ldots, \partial_{\Phi_0}^{\varepsilon_{2 g}}$ of the tangent space at the origin of $A \times A$. Then, the tangent space at the origin of $i(A)$ has as its basis the operators

$$
\delta_i=\partial_{\Phi_0}^{\varepsilon_i+N \epsilon_{g+i}}, \ i=1, \ldots, g
$$

\noindent
and the operator $\delta_{(P, [N]\cdot P)}^{k, T_{i(A)}}$ is given by

$$
\left(\partial_{\Phi_{(P, [N] \cdot P)}}^{\epsilon_1+N \epsilon_{g+1}}\right)^{k_1} \circ \cdots \circ\left(\partial_{\Phi_{(P, [N] \cdot P)}}^{\epsilon_{g}+N \epsilon_{2g}}\right)^{k_{g}}=\left(\delta_1 \circ t_{(P, [N] \cdot P)}^*\right)^{k_1} \circ \cdots \circ\left(\delta_{g} \circ t_{(P, [N] \cdot P)}^*\right)^{k_{g}},
$$

\noindent
where $k=\left(k_1, \ldots, k_{g}\right) \in \N^{g}$. 

\medskip

\noindent
We want $F$ to vanish at origin to order at least $T_0$, for sufficiently large $T_0 \in \N$. Therefore, the system we are trying to solve is
\begin{equation}\label{system siegel}
\delta_{(0, 0)}^{k} F=0, \quad \forall \ k \in \N^{g}, \quad|k| \leq T_0  .
\end{equation}

\medskip

% \section{Siegel Lemma}
% We write the system \eqref{system siegel} as

% $$
% \delta^{k} F(0) = 0, \quad \forall \ k \in \N^{g}, \quad|k| \leq T_0  .
% $$

\noindent
Note that above system gives rise to a $T_w \times L_w$ matrix with
    $$
    L_w = \binom{L+g}{g}^2 \sim \frac{L^{2g}}{(g  !)^2} \ \ \ \ \ \text{  and } \ \ \ \ \ T_w = \sum_{I \leq T_0}\binom{I+g}{g} \sim \frac{T_0^{g}}{g!}.
    $$
\medskip

% We use the following lemma to bound the elements of above matrix.

% \medskip
%      \begin{theorem}(\cite[Theorem 4.1]{David}) There exist constants $c_{17}, c_{18}$ and $c_{19}$ depending only on $g$ such that for any real number $h, \ h \geqslant 1$, any rational integer $\delta, \ \delta \geqslant 1$, any abelian variety $A(\tau)$ defined over a number field $K$ of degree over $\Q \leqslant \delta$, of height $h(A(\tau)) \leqslant h$, any polynomial $P$ with coefficients in $K$ in the abelian functions of degree $L$ and height $h(P)$, and any differential monomial $\Delta$ of $\mathcal{D}[T]$, the function $\Delta(P)$ is expressed as a polynomial in the $f_1, \ldots, f_N$ of degree at most $L+T$ and height

% $$
% h(\Delta P) \leqslant h(P)+c_{17} T h+c_{18}(L+T) \log (L+T) .
% $$
% \noindent
% Furthermore, the coefficients of $\Delta(P)$ are found in an extension $K(P)$ of $K$, of degree over $K \leqslant c_{19}$.
%  \end{theorem}
%  \begin{proof}
%      See \cite[Theorem 4.1]{David}.
%  \end{proof}
\medskip

\noindent
 Let $F=\sum_{\bf i} b_{\bf i} \mathbf{X}^{\bf i}$ be a polynomial with coefficients in $\bar{K}$. We define its height $h(F)$ as the logarithmic Weil height of the projective point defined by its coefficients and $1$. We need the following proposition to bound the elements of the coefficient matrix in \eqref{system siegel}.

%      \begin{proposition}[David \cite{David}]\label{David} There exist positive constants $c_{11}$ and $c_{12}$ depending only on $g$ such that for any real number $h, \ h \geqslant 1$, any abelian variety $A/K$ of height $h(A) \leqslant h$, any polynomial $P$ with coefficients in $K$ in the abelian functions of degree $L$ and height $h(P)$, and any differential monomial $\Delta$ of $\mathcal{D}[\textbf{T}]$, the function $\Delta P$ is expressed as a polynomial in the $f_1, \ldots, f_N$ of degree at most $L+T$ and height
% $$
% h(\Delta P) \leq h(P)+c_{11} T h+c_{12}(L+T) \log (L+T) .
% $$
% \noindent
% Furthermore, the coefficients of $\Delta P$ lie in a finite extension $K'$ of $K$, of degree over $K$ at most  $c_{13}$.
%  \end{proposition}
\begin{proposition}[David \cite{David}]\label{David}
Let $h\geq 1$, and let $A/K$
be an abelian variety of dimension $g$ with projective theta height \( h(A)\leq h
\). There exist positive constants $c_{11}$, $c_{12}$, and $c_{13}$, depending
only on $g$, such that the following holds. Let $P$ be a polynomial in the abelian functions
$f_1,\ldots,f_N$, with coefficients in $K$, degree $L$, and height $h(P)$. 
For every differential monomial $\Delta\in\mathcal{D}[\mathbf{T}]$ of
order at most $T$, the function $\Delta P$ can be expressed as a polynomial
in $f_1,\ldots,f_N$ of degree at most $L+T$ and satisfying
\[
h(\Delta P)
\leq
h(P)+c_{11}Th+c_{12}(L+T)\log(L+T).
\]
Moreover, the coefficients of this polynomial lie in a finite extension
$K'/K$ satisfying
\[
[K':K]\leq c_{13}.
\]
\end{proposition}
 \begin{proof}
     See \cite[Theorem 4.1]{David}.
 \end{proof}

\begin{lemma}\label{polynomial_bound}
 The height of each coefficient of system \eqref{system siegel} is bounded above by
$$
c_{21}\left(T_{0} \log \left(T_{0}+L\right)+T_{0} \log N+L\right)
$$
for an absolute constant $c_{21} > 0$ depending only on $A$ and $\rL$.
   
\end{lemma} 
\begin{proof} Recall that the operators $\delta_{(Q, N Q)}^{k}$ are defined as $\left(\delta_{1} \circ t_{(Q, N Q)}^{*}\right)^{k_{1}} \circ \cdots \circ\left(\delta_{g} \circ t_{(Q, N Q)}^{*}\right)^{k_{g}}$. Now from Proposition \ref{David}, we deduce that if $\delta$ is a monomial of weight $k$ in the $\delta_{\Phi_{0}}^{\epsilon_{j}}$ for $1 \leq j \leq 2 g$ and $F \in \Gamma\left(A \times A, \mathcal{M}_{\mid A \times A}^{\otimes L}\right)$, then $\delta(F)$ is a polynomial of degree $\leq L+k$ and

$$
h(\delta(F)) \leq h(F)+c_{12}(k+k \log (k+L)+k \log (N)).
$$

\medskip
\noindent
 Applying this to the operator $\delta_{1} \circ \cdots \circ \delta_{g}$ we obtain
 
$$
h\left(\delta_{1} \circ \cdots \circ \delta_{g}(F)\right) \leq h(F)+c_{13}(k \log (k+L)+k \log (N)) .
$$

\medskip
\noindent
Next, we use the fact that translations can be represented by forms of bi-degree $(2,2)$ to deduce that $t_{(P, [N] \cdot P)}^{*}(F)$ has degree $\leq 2 L$ and height $\leq h(F)+c_{14} L h(P, [N] \cdot P)$. Putting these two inequalities together and using the properties of the canonical height, we find that for $k \in \mathbb{N}^{g},|k| \leq T$ we have

$$
h\left(\delta_{(P, [N] \cdot P)}^{k}(F)\right) \leq c_{15}\left(h(F)+L+L N^{2} \hat{h}(P)+T_{0} \log \left(T_{0}+L\right)+T_{0} \log N\right) .
$$

\medskip
\noindent
Finally, we put $P=\mathcal{O}$ to deduce the desired result.
\end{proof}

\noindent
The following version of Siegel's lemma for a number field of degree $K$ and discriminant $D_{K/\Q}$ is due to Bombieri and Vaaler \cite[Theorem 9]{BV}.

\begin{lemma}(Siegel's Lemma)
\label{Siegel}
    Let $K$ be a number field of degree $d$ and discriminant $D_{K/ \Q}$. Let $A$ be an $M \times N$ matrix with entries in $K$ and $N > M$. Then there exist $N-M$ many linearly independent vectors ${\bf x}_l \in \mathcal{O}_K^N$ such that
    $$A\cdot{\bf x}_l = 0, \ \ \ \ \ \ l=1, \ 2, \ \ldots \ , \ N-M$$
    and
    $$
    \prod_{l=1}^{N-M} H({\bf x}_l) \leq |D_{K/\Q}|^{\frac{N-M}{2d}} H_{Ar}(A),
    $$
\end{lemma}
\noindent
where $H_{\text{Ar}}(A)$ is the arakelov height of $A$.

\medskip
\noindent
Note that in the lemma above, $H(\cdot)$ and $H_{\text{Ar}}(\cdot)$ denote the multiplicative Weil height and the multiplicative Arakelov height, respectively, defined by $H(\cdot) := e^{h(\cdot)}$ and $H_{\text{Ar}}(\cdot) := e^{h_{\text{Ar}}(\cdot)}$. Now we are ready to prove the existence of a non-zero solution of small height to the system \eqref{system siegel}. 

\begin{proposition}\label{Size of polynomial}
    If $L_w > T_w$, then the system \eqref{system siegel} has a non-zero solution in $\mathcal{O}_{K}^{L_w}$ such that
    $$
    h(F) \leq C \frac{T_w}{L_w - T_w} \left(T_0 \log \left(T_0+L\right)+T_0 \log N+L\right)
    $$
    for some $C > 0$ depending only on $A$ and $K$.
\end{proposition}
\begin{proof}
    The proof is an easy consequence of Lemma \ref{polynomial_bound} and Lemma \ref{Siegel}. Denote by $A$ the $L_w \times T_w$ matrix corresponding to the homogeneous system \eqref{system siegel} of linear equations and let $r$ be its rank. Consider a sub-matrix $B$ of a consisting of $r$ linearly independent rows of $A$. It follows from the definition of arakelov height that $h_{\text{Ar}}(A) = h_{\text{Ar}}(B)$. 
    % \text{col}(A)$ be the column space of $A$. Let $B$ be a submatrix Clearly, $\text{dim}_{\overline{\Q}} \left(\text{null}(A)\right) = L_w - \text{rank}_{\overline{\Q}}(A) \geq L_w - T_w$.

    \medskip
    
    Let $\textbf{y}_1, \textbf{y}_2, \ldots, \textbf{y}_r$ be the row vectors of $B$. Using Lemma \ref{polynomial_bound}, we get
    \begin{align*}
    h_{\text{Ar}}(A) = h_{\text{Ar}}(B) \leq \sum_{i=1}^r h_{\text{Ar}}(\textbf{y}_i) & \leq r \cdot \operatorname{max}_i h_{\text{Ar}}(\textbf{y}_i)\\
        & \leq c_{21} r \left(T_{0} \log \left(T_{0}+L\right)+T_{0} \log N+L\right). 
    \end{align*}

\noindent
    This along with Lemma \ref{Siegel} implies the existence of a non-zero solution $\textbf{x}$ of \eqref{system siegel} such that 
\begin{align*}
    h(\textbf{x}) &\leq \frac{\log |D_{K/\Q}|}{2d} + c_{21} \frac{r}{L_w - r} \left(T_{0} \log \left(T_{0}+L\right)+T_{0} \log N+L\right)\\   
    &\leq C \frac{ r}{{L_w - r}} \left(T_{0} \log \left(T_{0}+L\right)+T_{0} \log N+L\right).
\end{align*}
Now using the trivial upper bound on the rank $r$, i.e. the number of equations $T_w$, we get the desired solution. 
% From here using a remark by Roy-Thunder
% (\cite[p. 7]{Roy}) or a result of Bombieri-Vaaler \cite[Theorem 9]{BV} one can deduce integral solution such that $h(F)$ is
% bounded by the same quantity up to a constant.
\end{proof}

\section{\bf Height estimate}

Throughout this section, we fix a N\'{e}ron model of the abelian variety $A$ over $\mathcal{O}_K$. We also fix a basis of algebraic derivations on $A$. For all prime ideals $\mathfrak{p}$ of $\mathcal{O}_K$ larger than a constant depending only on $A$, this basis of derivations remains a basis of derivations modulo $\mathfrak{p}$. We have fixed a basis $(s_i)$ of $\Gamma^0(A, L)$ and we set $f_i := s_i/s_0$, which are affine functions on the open subset of $A$ where $s_0$ does not vanish.

\begin{proposition}
    There exists a basis of derivations $(\partial_{1}, \ldots, \partial_{g})$ on $A$ such that
$$
\forall (i,j): \quad \partial_{j} f_{i}
    = \sum_{(k,l)} y_{k,l}^{i,j} \, f_{k} f_{l}
$$
where the $y_{k,l}^{i,j} \in \overline{\mathbb{Q}}$ have height bounded only in terms of $A$.
\end{proposition} 
\begin{proof}
    See \cite[Theorem 2.9]{Galateau}.
\end{proof}
\medskip

\noindent
Above proposition implies that for all but finitely many places $v \in M_K$, we have 
\begin{equation}\label{integral}
   |y_{k,l}^{i,j}|_v \leq 1 
\end{equation}
for all $i,j,k$ and $l$. Moreover, we know that there exists a finite set of rational primes $\mathcal{P}$ such that $A/K$ has good reduction at every place above the rational primes $p$ not contained in $\mathcal{P}$. By enlarging the set $\mathcal{P}$ if necessary, we can further assume that for any prime $p$ outside $\mathcal{P}$, the abelian variety $A/K$ has good reduction at every place $v$ above $p$, and each such $v$ satisfies condition \eqref{integral}.
\medskip
% The isomorphism:
% $$
% \hat{\mathcal{A}}_{\mathcal{O}_v} \simeq \mathcal{O}_v [[ T_1, \ldots, T_g ]],
% $$

% \noindent
% gives the embedding
% $$
% H^0\left(\mathcal{A}_{\mathcal{O}_v}, \mathcal{O}_{\mathcal{A}_{\mathcal{O}_v}}\right) \rightarrow \mathcal{O}_v [[ T_1, \ldots, T_g ]]. $$

\noindent
The section $F$ constructed using the system \eqref{system siegel} has a Taylor expansion around $0$ and is given by 

% of $H^0\left(\mathcal{A}_{\mathcal{O}_v}, \mathcal{O}_{\mathcal{A}_{\mathcal{O}_v}}\right)$ 

$$
F(T_1, T_2, \ldots, T_g) = \sum_{\mu \in \mathbb{N}^g} \partial^\mu F (0) \textbf{T}^\mu,
$$
\noindent
where $\mu = (\mu_1, \mu_2, \ldots, \mu_g)$ and $\textbf{T}^\mu = T_1^{\mu_1}T_2^{\mu_2} \ldots T_g^{\mu_g}$. By the construction, it follows that 
$\partial^\mu F \left( 0 \right) = 0
$ for all $\mu$ such that $|\mu| \leq T_0.$ 
\medskip

\noindent
Let $P \in A(\overline{K})$ and let $w$ be a place of $K(P)$ lying above a rational prime $p \notin \mathcal{P}$. Suppose that $P$ reduces to the identity component under reduction modulo $w$ in the integral model. This implies
$$
\left|T_i \left(P\right)\right|_w \leq q^{-1} \ \  \forall \ \ 1 \leq i \leq g,
$$
where $q$ is norm of $w$.

\medskip

\noindent
Moreover, for all $\mu \in \mathbb{N}^g$, using \cite[Section 5B]{Galateau}, we have
$$
|\partial^\mu F(0)|_w \leq 1
$$

% \medskip
% \noindent
% Therefore, the Taylor expansion of $F$ can be written as
% $$ F(\mathbf{T}) = \sum_{\mu \in \N^g, |\mu| > T_0} \partial^\mu F(0)  \ \mathbf{T}^\mu. $$

\medskip
\noindent
Let $T_f$ be the order of smallest derivative $\Delta_f$ such that $\Delta_f F(P) \neq 0$. If $P \equiv O  \ (\operatorname{mod} {w})$, this gives

$$ | \Delta_f F(\mathbf{T})(P)|_w = \Bigg|\sum_{\mu \in \N^g, |\mu| > T_0} \partial^\mu F (0) \ \Delta_f \mathbf{T}^k (P)\Bigg|_w \leq q^{-(T_0 - T_f) } 
$$
or
\begin{equation}\label{upper}
 \log | \Delta_f F(\mathbf{T})(P)|_w \leq -(T_0 - T_f) \ \log q.
\end{equation} 
\medskip

\noindent
On the other hand, using Proposition \ref{David}, we have 
$$
h(\Delta_f F(P)) \leq c_{11}(h(F) + L(1 + N^2) \hat{h}(P) + T_f \log(T_f + L) + T_f \log N),
$$
where $c_{11}$ depends only on $A$ and $K$.
\medskip

\noindent
We now apply Liouville's inequality \eqref{Liouville's inequality} to obtain
\begin{equation}\label{lower}
 \sum_{w|v,\ w \in M_{K(P)}} \log |\Delta_f F(P)|_w \geq - c_{11}[K(P) : K] (h(F) + L(1 + N^2) \hat{h}(P) + T_f \log(T_f + L) + T_f \log N).
\end{equation}

\noindent
 In \eqref{upper}, considering the contribution of all places of $K(P)$ of norm $q$ and using \eqref{lower}, we obtain
\begin{equation}\label{lower-2}
    c_{11}(h(F) + L(1 + N^2) \hat{h}(P) + T_f \log(T_f + L) + T_f \log N) \geq (T_0-T_f)  \frac{\tilde{N}_q(K(P))}{[K(P) : \Q]}  \log q,
\end{equation}
where $\tilde{N}_q(K(P))$ is the number of places of norm $q$ in $K(P)$ such that $P \equiv O \ ( \text{mod} \ w)$.

\medskip
\noindent
Recall that from Proposition \ref{Size of polynomial}, we have  
    $$
    h(F) \leq C \frac{T_w}{L_w - T_w} \left(T_0 \log \left(T_0+L\right)+T_0 \log N+L\right) ,
    $$
where $T_w \sim {T_{0}}^{g}/g!$ and $L_w \sim L^{2g}/(g!)^2$. Therefore, plugging this in \eqref{lower-2}, we obtain 

\begin{align}\label{final1}
    (T_0-T_f)  \frac{\tilde{N}_q(K(P))}{[K(P) : \Q]}  \log q &\leq C_1 \Bigl( \frac{T_w}{L_w - T_w} \left(T_0 \log \left(T_0+L\right)+T_0 \log N+L \right) + L(1 + N^2) \hat{h}(P)  \nonumber \\ &+ T_f \log(T_f + L) + T_f \log N  \Bigl).
\end{align}

\subsection{\bf Setting the parameters}
For points $P \in A(\mathcal{K})$, where $\mathcal{K}$ is an asymptotically positive extension with $\psi_q > 0$, the quantity $\tilde{N}_q(K(P)) / [K(P) : \mathbb{Q}] \gg \psi_q$ for sufficiently large $[K(P) : K]$. We set the parameters $L$, $N$, and $T_f$ as follows:
\begin{align*}
    L=\sqrt{T_0} \log T_0,\ N^2 = \frac{\sqrt{T_0}}{\log T_0}, \ T_f = \frac{T_0}{(\log T_0)^2}.
\end{align*}

\noindent
Hence, $T_0$ remains a free parameter, to be chosen sufficiently large. With these choices of parameters, we get
\begin{align*}
L_w \sim \frac{L^{2g}}{(g!)^2} = \frac{{T_0}^g}{(\log T_0)^{2g}(g!)^2} \sim T_w \frac{(\log T_0)^{2g}}{g!} \gg T_w.
\end{align*}

\noindent
Therefore, the hypothesis of Proposition \ref{Size of polynomial} holds for a sufficiently large value of $T_0$. Moreover, for large values of $T_0$, we have
$$
\frac{T_w}{L_w - T_w} \ll_g \frac{1}{(\log T_0)^2}, \ \ \ \ \ \ \ \ L(1+N^2) \sim T_0, \ \ \ \ \ \ \ \  \log (T_0 + L) \ll \log T_0
$$

$$
\log (T_f + L) \ll \log T_0,  \ \ \ \ \ \ \ \ \log N \ll \log T_0.
$$

% \begin{enumerate}[(a)]
%     \item $$
%     \rho^g \ T_0 \log L \leq T_0  \psi_q \log q \implies \rho^g \log L \leq \psi_q \log q,
%     $$
%     \item $$
%     \rho^g L \leq T_0 \psi_q \log q \implies  \frac{\rho^{g-1}}{L} \leq \psi_q \log q,
%     $$
%     \noindent
%     and
    
%     \item $$
%     T_f \log L \leq  (\psi_q \log q) \ T_0 \implies T_f \leq \operatorname{min}\left( \frac{T_0}{2}, \rho \frac{L^2}{\log L} \psi_q \log q \right).
%     $$ 
% \end{enumerate}

\subsection{\bf Lower bound on height}\label{subsection}

The choices of parameters altogether imply that if $P \equiv O$ for at least $\psi_q \cdot [K(P) : \mathbb{Q}]$ many places of norm $q$, and if the order of vanishing of $F$ at $P$ is bounded by $T_f$, then for a sufficiently large value of $T_0$, equation \eqref{final1} yields
\begin{equation}\label{height bound}
   \hat{h}(P) \geq C_2 \cdot \psi_q \log q, 
\end{equation}
where $C_2$ depends only on $A$, $K$ and $T_0$.

\section{ \bf Zero estimate}\label{Sec zero}

Throughout this section $A/K$ denote a geometrically simple abelian variety defined over the number field $K$ and $q$ denote a prime power $p^f$ such that $p \notin \mathcal{P}$.  For a fixed finite extension $E$ of $K$, define the set
$$
% Z = \{ P \in A(L) \ | \ P \equiv O \ (\text{mod}\ w) \  \forall \ w |v \} \quad \text{and} \quad 
Z_F = \{ P \in A(E) \ | \ \Delta F(P) = 0 \ \text{for} \ |\Delta| \leq T_f  \}.
$$

\noindent
Let $P \in A(E)$ such that $\hat{h}(P) < C_2 \cdot \psi_q \log q$. Our height estimate \eqref{height bound} implies that if $P \equiv O$ for at least $\psi_q [E:\Q]$ many places of $E$ of norm $q$, then $P \in Z_F$. \medskip

\noindent
We need the following result due to P. Philippon  for the zero estimate.

\begin{theorem}[{\cite[Philippon]{Philippon}}]\label{Philippon}
    
Let $A$ be an abelian variety defined over $\C$, of dimension $g$, embedded in a projective space $\mathbb{P}^n$. Denote by deg($\cdot$) the degree map on the closed subsets of ${\mathbb{P}}^n$. Let $V$ be a subspace of $T_A(\C)$, the tangent space at the origin of $A$, and $S$ be a finite set of points of $A$. We denote by $\Gamma^{(g)}(S)$ the set: $\left\{\Sigma_{i=1}^g P_i, P_i \in S\right\}$. Consider, furthermore, an element $F$ of $\C\left[X_0, \ldots, X_N\right]$, homogeneous of degree $D$, which vanishes to order $gT+1$ along $V$ on $\Gamma^{(g)}(S)$, but which is not identically zero on A. Then there exists an abelian subvariety $B$ of $A$, distinct from $A$ defined by equations of degree at most $c_0 D$, such that:

$$
\begin{aligned}
& \binom{T+\operatorname{codim}_{V}\left(V \cap T_B\right)}{\operatorname{codim}_{V}\left(V \cap T_B\right)} \cdot \operatorname{card}\left(\frac{S+B}{B}\right) \cdot \operatorname{deg}(B) D^{\operatorname{dim}(B)}  \leq \operatorname{deg}(A) (2 D)^{\operatorname{dim}(A)}.
\end{aligned}
$$
\end{theorem} 

\bigskip

% \noindent
% Assume that $\Delta_f F (P) \neq 0$ for all $P \in {\mathcal{S}}^{(i)} \setminus S_F(E)$, hence for such points we get 

% $$
% \hat{h}(P) \ll \rho \log q.
% $$
\noindent
Consider the set 
\begin{equation}\label{Vanishing set}
  S = \Big\{ P \in A(E) \ | \ P \equiv O \ \text{for at least} \ \psi_q [E:\Q] \ \text{many places of norm $q$} \ \ \& \ \ \ \hat{h}(P) < C_2 \cdot \psi_q \log q \Big\}.  
\end{equation}

\medskip

\noindent
The function $F$ is a polynomial of degree at most $L$ in the abelian functions of the product variety $A \times A$. Moreover, $F$ vanishes on all the points of the form $(Q, N Q)$ to order $T_f$ along the subvariety $i(A)$ of $A \times A$ defined by the relation $z_2=N z_1$. Therefore, Theorem \ref{Philippon} ensures the existence of a proper algebraic subgroup $B$ of $A$ such that the following inequality holds:

% Since $\mathcal{K}$ is a totally $v$-adic extension of $K$.

% Using the formulas for multiplications by $N$, to a function $G$ on $A$, of degree $\leqslant c' L N^2$, which vanishes on $Z$, to order $T_0$. Moreover, it is not identically zero on $A$.

% : if this were the case, we would deduce that the function $F$ vanishes identically on $H$ (therefore along $H$ at any order). We can then apply the zero lemma (Theorem 6.4) with any value of $T$. There therefore exists a proper subgroup $B$ of $A \times A$ such that:

% $$
% \begin{aligned}
% & \binom{T+\operatorname{codim}_{T_H}\left(T_B \cap T_H\right)}{\operatorname{codim}_{T_H}\left(T_B \cap T_H\right)} \cdot \operatorname{deg}(B) L^{\operatorname{dim}(B)} \leq c_{41} \operatorname{deg}(A) \cdot(2 L)^{2 \operatorname{dim}(A)}.
% \end{aligned}
% $$

% \medskip
% \noindent
% For this inequality to be satisfied $(\forall T)$, $B$ must contain $H$. Since $A$ is simple, $B=H$. 
% \medskip
% \noindent
% We deduce:
% $$
% \operatorname{deg}(H) \leqslant c' L^{2 g},
% $$
% \noindent
% that is,

% $$
% N^2 \leqslant c'' L,
% $$
% which is incompatible with the choice of parameters.

% \medskip

\begin{align}\label{zero estimate}
\binom{T_0+\operatorname{codim}\left(B\right)}{\operatorname{codim}\left(B\right)} &\cdot \operatorname{card}\left( \frac{Z + B}{B} \right) \cdot \operatorname{deg}(B) L^{\operatorname{dim}(B)} \nonumber \\ 
& \leq \operatorname{deg}(A)\cdot(2LN^2)^{\operatorname{dim(A)}}.
\end{align}
Since $A$ is simple, the only proper subgroup of $A$ is $\{O\}$, thus $B$ is trivial. This gives $\textrm{codim}(B)= g$ and $\textrm{deg}(B)=0$. Hence
$$
{T_f}^g \cdot \#S \ \leq C_3 \ L^{2g},
$$
where $C_3$ depends only on $A$ and $K$.
\medskip

\noindent
We deduce by replacing different parameters in above inequality with their values to obtain:
$$
 \# S \ \ \leq \ C_3\ \frac{L^{2g}}{T_f^g}  \leq C_3 (\log T_0)^{4g}.
$$

\noindent
Hence, we have the following proposition.

\begin{proposition}\label{Bogomolov}
    Let $A/K$ be a geometrically simple abelian variety and $E/K$ be a finite extension of $K$. Then, the set $S$ defined in \eqref{Vanishing set}
    % $Z = \{P \in A(\mathcal{K}) \, | \ P \equiv O \, (\operatorname{mod} \, $w$) \, \, \forall \, w \in M_{K(P)} \, \, \textrm{above} \, \, v \}$
    % such that for $P \in {\mathcal{S}}^{(i)} \setminus S_F(E)$ 
    % $$
    % \hat{h}(P) \ll \rho \log q
    % $$
    has cardinality at most
    $$
     C_3 \cdot (\log T_0)^{4g},
    $$
    where $C_3$ is a constant depending only on $A$ and $K$.
\end{proposition}

% \section{Bogomolov property for totally $v$-adic extension}

% Let $A/K$, $\mathcal{K}$ and $v$ be as in Proposition \ref{Bogomolov}. We claim that the set 
% $$
% S = \{P \in A(\mathcal{K} \, \, | \, \, \hat{h}(P) \ll \rho \frac{\log q}{q^{2g}}\}
% $$
% is finite.

% Let $P \in S$ and $w$ be any place of $K(P)$ above $v$. The reduction of abelian variety at $w$ has the cardinality $\# A(\mathbb{F}_q) \sim q^{g}$. Therefore, 
% $$ \#A(\mathbb{F}_q) \cdot P \equiv \mathcal{O} \, \, (\operatorname{mod} w)$$ 
% for all $w \mid v$.

% \medskip
% \noindent
% Hence the set
% $$
% S' = \#A(\mathbb{F}_q) \cdot S = \{\#A(\mathbb{F}_q) \cdot P \, | \, P \in S\}
% $$
% satisfy the property that $Q \equiv 0$ for all places $w\mid v $ of $K(Q)$ and $\hat{h}(Q) = \hat{h}(\#A(\mathbb{F}_q) \cdot P) \ll \rho \log q$.

% \medskip

% Thus from the Proposition \ref{Bogomolov}, 

% $$
% |S'|  \ll \frac{1}{\rho^{g(g+1)}}.
% $$
% We conclude from here that 
% $$
% |S| \ll \frac{q^{2g}}{\rho^{g(g+1)}}. 
% $$

% \begin{theorem}
%     Let $A/K$ be an abelian variety and $\mathcal{K}/K$ be a totally $v$-adic extension of $K$. Then, the set
%     $$
%     S = \{P \in A(\mathcal{K} \, \, | \, \, \hat{h}(P) \ll \rho \frac{\log q}{q^{2g}} \}
%     $$
%     has cardinality 
%     $$
%     \ll \frac{q^{2g}}{\rho^{g(g+1)}},
%     $$
%     where $q$ is norm of $v$ and the implied constants depend only on $A$, $\mathcal{L}$ and  $K$.
% \end{theorem}

\section{\bf Proof of Theorem \ref{Main 2}}

As in the previous sections, unless otherwise stated, $q$ denotes a prime power $p^f$ with $p\notin\mathcal{P}$. We first prove Corollary \ref{Bogomolov_totally_p-adic} independently and then extend these ideas to prove Theorem \ref{Main 2}. The main difficulty in extending the argument from Corollary \ref{Bogomolov_totally_p-adic} to Theorem \ref{Main 2} is the need for a suitable box principle that allows us to combine the contributions from different places linearly.

\begin{proof}[Proof of Corollary \ref{Bogomolov_totally_p-adic}]

Let $A/K$ be a geometrically simple abelian variety and $\mathcal{K}/K$ be an asymptotically positive extension with $\psi_q(\mathcal{K}) > 0$. For a non-negative real number $\epsilon > 0$, define the set

$$
S = \{P \in A(\mathcal{K}) \,\, | \,\, \hat{h}(P) < \epsilon\}.
$$

\medskip
\noindent
It is enough to show that the set $S$ is finite for some $\epsilon > 0$. Assume, contrary to this, that the set described above is infinite. We take a large finite extension $E$ of $K$ contained in $\mathcal{K}$. For any place $w \in M_E$ of norm $q$, consider the maps

$$
A/\mathcal{O}_E \xrightarrow[]{[\# A(\F_q)]} A/\mathcal{O}_E \xrightarrow[]{\text{mod} \ w} A/\F_q,
$$

\medskip
where the first map is multiplication by $\# A(\F_q)$ and the later map is reduction modulo $w$. Clearly, every point $P \in A$ maps to the identity of $A(\mathbb{F}_q)$ under this composition. Since the multiplication map is surjective and has degree $\#A(\mathbb{F}_q)^{2g}$, we conclude that the resulting set

$$
S' = \{ [\# A(\F_q)] \cdot P \ | \ P \in S   \}
$$

\medskip

\noindent
is infinite. Moreover, every $P \in S'$ satisfies $\hat{h}(P) < q^{2g} \cdot \epsilon$, using the Hasse bound (i.e., $\# A(\mathbb{F}_q) \sim q^g$) and the fact that $\hat{h}$ is a quadratic form on $A(\overline{K})$.
% and $P \in S'$ satisfy $\hat{h}(P) < q^{2g} \cdot \epsilon$ using the Hasse's bound (i.e. $\# A(\F_q) \sim q^g$) and the fact that $\hat{h}$ is a quadratic form on $A (\overline{K})$.
\medskip

It follows from the discussion above that every point $P \in S' \cap A(E)$ reduces to $O$ modulo $w$ for all places $w$ of $E$ of norm $q$.  
Let $E$ be as above. Since $\mathcal{K}$ is asymptotically positive, for a sufficiently large $[E : K]$ we have
\[
\frac{N_q(E)}{[E : \mathbb{Q}]} \gg \psi_q .
\]

Choosing $\epsilon = C_2 \cdot \psi_q \frac{\log q}{q^{2g}}$ and applying Proposition \ref{Bogomolov}, we obtain
\[
|S' \cap A(E)| \le C_3 (\log T_0)^{4g}.
\]

This bound holds for arbitrarily large $E$, so $S'$ must be finite. Because multiplication by $\# A(\mathbb{F}_q)$ is a finite morphism, the set $S$ is also finite. Consequently, $A(\mathcal{K})$ has the Bogomolov property. 
\end{proof}
\medskip

\noindent
Now our next goal is to prove Theorem \ref{Main 2}. The strategy is similar to the proof of Corollary \ref{Bogomolov_totally_p-adic} but a more careful zero estimate and a box principle is require to obtain desire result.  Fix a real number $X > 0$. Note that throughout the proof $\sum'_q$ and $\prod'_q$ denote sum and product over prime power $q=p^f$ such that $p \notin \mathcal{P}$. For a finite extension $E/K$ and a point $P \in A(E)$, define
$$
\phi_X(P, A(E)) = \sum'_{\substack{w \in M_E, \ q_w\leq X \\ P \equiv O \ (\text{mod} \ w)}} \log q_w,
$$
where $q_w$ denote the norm of $w$ and $\sum'$ denote sum over those $q_w$ which are power of prime lying outside $\mathcal{P}$.
\medskip
\noindent
Using this in \eqref{upper}, we get

\begin{equation}\label{general1}
 \sum'_{\substack{w \in M_E, \ q_w\leq X }} \log | \Delta_f F(\mathbf{T})(P)|_w \leq -(T_0 - T_f) \ \phi_X(P, A(E))
\end{equation} 

\noindent
Now using Liouville's inequality \eqref{Liouville's inequality} in \eqref{general1}, we obtain
\begin{equation}\label{general-2}
    c_{11}(h(F) + L(1 + N^2) \hat{h}(P) + T_f \log(T_f + L) + T_f \log N) \geq (T_0-T_f)  \frac{\phi_X(P, A(E))}{[E : \Q]}  \log q.
\end{equation}

\subsection{\bf Setting the parameters}
Fix a real number $X > 0$. As earlier, we set the parameters $L$ and $T_f$ as follows:
\begin{align*}
    L=\sqrt{T_0} \log T_0,
    % T_0 = \frac{L^2}{(\log T_0)^2},
    \ T_f = \frac{T_0}{(\log T_0)^2}.
\end{align*}

\noindent
We slightly modify the parameter $N$ as 
    $$
     \ N^2 = \frac{1}{\left( \prod_{q \leq X}'  q \right)^{2g}}  \cdot \frac{\sqrt{T_0}}{\log T_0},
    $$
where the product runs over all prime power $q = p^r$ such that $p \notin \mathcal{P}$. 

% we need $N^2 \sim L$, $\log T_f =  O(\log T_0) = O(\log L)$ and $(1 -T_f/T_0) \sim \frac{1}{2}$. Let $T_0 = \rho L^2$ and choose the parameters as follows

% \begin{enumerate}[(a)]
%     \item $$
%     \rho^g \ T_0 \log L \leq T_0  \sum'_{q \leq X} \psi_q \log q \implies \rho^g \log L \leq \sum'_{q \leq X} \psi_q \log q,
%     $$
%     \item $$
%     \rho^g L \leq T_0 \sum'_{q \leq X} \psi_q \log q \implies  \frac{\rho^{g-1}}{L} \leq \sum'_{q \leq X} \psi_q \log q,
%     $$
%     \noindent
%     and
    
%     \item $$
%     T_f \log L \leq T_0 \sum'_{q \leq X} \psi_q \log q \implies T_f \leq \operatorname{min}\left( \frac{T_0}{2}, \rho \frac{L^2}{\log L} \sum'_{q \leq X} \psi_q \log q \right).
%     $$ 
% \end{enumerate}

\medskip

\subsection{\bf Height estimate}
Similar to Section \ref{subsection}, the choices of parameters altogether imply that if $\phi_X(P, A(E))\gg \left( \sum'_{q \leq X} \psi_q \log q \right)[E : \Q]$, and the order of vanishing of $F$ at $P$ is bounded by $T_f$, then for a sufficiently large value of $T_0$, equation \eqref{general-2} yields 
$$
\hat{h}(P) \geq C_2 \cdot \left( \prod_{q \leq X}'  q \right)^{2g} \cdot \sum'_{q \leq X} \psi_q \log q,
$$
where $C_2$ depends only on $A$ and $K$.

\subsection{\bf Zero estimate}\label{Zero estimate_2}
% For $E/K$ as above, using a zero estimate similar to the earlier case, we conclude that all the points in the set 
% $$
% Z = \Big\{ P \in A(E) \ | \ \frac{\phi_X(P, A(E))}{[E : \Q]} \gg \sum'_{q \leq X} \psi_q \log q \Big\} \quad \text{and} \quad Z_F = \{ P \in A(E) \ | \ \Delta F(P) = 0 \ \text{for} \ |\Delta| \leq T_f  \}.
% $$

% \noindent
% This allow us to conclude that $P \in Z$ satisfying $\hat{h}(P) \leq C_2 \cdot \sum'_{q \leq X} \psi_q \log q$ lie in $Z_F$. Now following the similar zero estimate as in Section \ref{Sec zero}, we get  
% 
For $E/K$ as above, define the set
$$
S = 
\Bigg\{
P \in A(E) \ | \ \frac{\phi_X(P, A(E))}{[E : K]} \gg \sum'_{q \leq X} \psi_q \log q \ \ \& \ \ \ \hat{h}(P) < C_2 \cdot \left( \prod_{q \leq X}'  q \right)^{2g} \cdot \sum'_{q \leq X} \psi_q \log q 
\Bigg\}.
$$
A zero estimate similar to the earlier case implies that every element of $S$ vanishes at $F$ to order at least $T_f$. Moreover, cardinality of $S$ is bounded above by
\begin{equation}\label{zero_2}
     C_3 \cdot (\log T_0)^{4g},
\end{equation}
where $C_3$ depends only on $A$ and $K$.

\begin{proof}[Proof of Theorem \ref{Main 2}]
    
Let $\rK/K$ be an infinite extension and $\{K_i\}_{i \geq 1}$ be a tower of number fields such that $\mathcal{K}= \cup_{i \geq 0} K_i$. For any real number $\epsilon >0$ such that  $\epsilon \leq \frac{C_2}{2} \cdot \sum'_{q \leq X} \psi_q \log q$, we consider the set 
$$
\Bigg\{ P \in A(\mathcal{K}) \ | \ \hat{h}(P) < \epsilon \Bigg \}.
$$
\noindent
Assume that above set is infinite and we denote by $\tilde{S}^{(i)}$ its $K_i$ rational points defined as
$$
\tilde{S}^{(i)}= \Bigg\{ P \in A(E) \ | \ \hat{h}(P) < \epsilon \Bigg \}.
$$
% \begin{proposition}\label{Main_2 key}
%     If $S(L) \rightarrow \infty $ as $[L : \Q] \rightarrow \infty$ then the set 
%     \begin{equation}\label{condition}
%          \Bigg\{ P \in Z \ \Bigg\mid \ \frac{\phi_X(P, A(L))}{[L : K]} \gg \sum_{q \leq X} \psi_q \log q \Bigg\}
%     \end{equation}
%     has cardinality at least $(\prod_{q \leq X} q)^g \cdot |S'(L)|$. 
% \end{proposition}
% \begin{proof}
%     It is clear that $\#Z \geq \# S(L)$. We know that $\#S(L) \rightarrow \infty$ as $L$ is taken large. Consider 
%     \begin{align}
%     \sum_{P \in Z} \quad \sum_{\substack{w \in M_E, \ q_w \leq X \\ P \equiv O \ (\text{mod} \ w)}}   \frac{\log q_w}{[L:\Q]}  &\geq  \sum_{w \in M_E, \ q_w \leq X}  \sum_{\substack{P \in Z, \\ P \equiv O \ (\text{mod} \ w)}} \frac{\log q_w}{[L:\Q]} \nonumber \\
%     &\gg \# S(L) \sum_{q \leq X} \frac{\psi_q}{q^{2g}} \log q
% \end{align}
% \end{proof}
% \begin{proof}
%     Consider the map 
%     $$
%     \overline\rho \ :\ A(L) \longrightarrow \prod_{w \in M_E, \ q_w \leq X} A(\F_{q_w}).
%     $$
%    Using the fact that RHS has cardinality $\sim (\prod_{q \leq X} q)^g$, at least $\#S'(L)/(\prod_{q \leq X} q)^g$ many points of $S'(L)$ lie in the fibre of a fix point under $\overline\rho$ and therefore $Z$ has a subset of cardinality at least $\#S'(L)/(\prod_{q \leq X} q)^g$ satisfying \eqref{condition}.
% \end{proof}
Consider the set ${\mathcal{S}}^{(i)}: = \tilde{S}^{(i)}-\tilde{S}^{(i)}$. Denote by ${\mathcal{S}}^{(i)}_F$ to be the set of points of ${\mathcal{S}}^{(i)}$ such that the order of vanishing of $F$ at $P$ is larger than $T_f$. Consider the set 
$$
{\mathcal{W}_F}^{(i)}: = \left(\prod_{q \leq X}' \#A(\F_q) \right) \cdot {\mathcal{S}}^{(i)}_F = \Big\{ \Big[ \prod_{q \leq X}' \#A(\F_q) \Big] \cdot P \  | \ P \in {\mathcal{S}}^{(i)}_F. \Big\} 
$$

Clearly, every element $P$ satisfy the congruence condition $P \equiv \mathcal{O} \ (\operatorname{mod} \ w)$ for every place $w \in M_E$ of norm $q \leq X$, where $q$ is a power of prime $p \notin \mathcal{P}$. Therefore, every element $P \in {\mathcal{W}}^{(i)}_F$ satisfy 
$$
\ \frac{\phi_X(P, A(E))}{[E : K]} \gg \sum'_{q \leq X} \psi_q \log q \ \ \text{and} \ \ \ \hat{h}(P) < C_2 \cdot \left( \prod_{q \leq X}'  q \right)^{2g} \cdot \sum'_{q \leq X} \psi_q \log q.
$$
Hence, using the zero estimate in Section \ref{Zero estimate_2}, we obtain
$$
|{\mathcal{W}_F}^{(i)}| \leq C_3 \cdot \left( \log T_0 \right)^{4g}.
$$

Using the fact that $\prod_{q \leq X}' \#A(\F_q) \sim \left( \prod'_{q \leq X} q \right)^{g}$ and multiplication by an integer $m$ is a finite endomorphism on abelian variety of degree $m^{2g}$, we deduce

\begin{equation}
|{\mathcal{S}}^{(i)}_F| \leq C_3 \cdot \left( \prod'_{q \leq X} q\right)^{2g^2} \cdot (\log T_0)^{4g}.  \nonumber
\end{equation}

% Assume that $S$ is infinite and consider a large number field $E/K$ lying inside $\mathcal{K}$. Denote by $S(E)$ the set of $E$-rational points of $A$ lying in $S$. It follows from Dirichlet's box principle that the number of zeros of $F$ to order $ > T_f$ in $S(E)$ can not be more than $C_3 \cdot \left( \prod'_{q \leq X} q\right)^{2g^2} \cdot (\log T_0)^{4g}$. 
We conclude from the above inequality that the number of zeros of $F$ to order at least $T_f$ inside the set ${\mathcal{S}}^{(i)}$ is bounded by a quantity independent of the field $K_i$. For the rest of points $P \in {\mathcal{S}}^{(i)} \setminus {\mathcal{S}}^{(i)}_F$, the product formula gives
$$
\sum_{w \in M_{K_i}} \log |\Delta_f F(P)|_w = 0.
$$

\noindent
Summing this equality over all points $P \in {\mathcal{S}}^{(i)} \setminus {\mathcal{S}}^{(i)}_F$ with multiplicity, we obtain 
$$
\sum_{P \in {\mathcal{S}}^{(i)} \setminus {\mathcal{S}}^{(i)}_F} \sum_{w \in M_{K_i}} \log |\Delta_f F(P)|_w = 0.
$$
% If the number of zeros of $F$ in $Z$ is $\gg \frac{1}{\rho^{g(g+1)}}$, we are done. Therefore assume that $F$ has at most $\frac{1}{\rho^{g(g+1)}}$ many zeros in $Z$, call it $Z_F$.

% Let $M_E^{\textrm{good}}$ denote the set of all non-archimedean places of $L$ at which $A$ has good reduction. 
\noindent
The Northcott theorem implies that ${\mathcal{S}}^{(i)}$ is finite and thus interchanging the summation and collecting only the desired places of norm $ q \leq X$, we can write

\begin{align*}
\sum_{P \in {\mathcal{S}}^{(i)} \setminus {\mathcal{S}}^{(i)}_F} & \sum_{w \in M_{K_i}} \log |\Delta_f F(P)|_w  \, =  \, \sum_{P \in {\mathcal{S}}^{(i)} \setminus {\mathcal{S}}^{(i)}_F} \left( \sum'_{\substack{w \in M_{K_i}, \\  q_w \leq X}} \log |\Delta_f F(P)|_w \, + \,  \sum_{\substack{w \in {M_{K_i}}^\infty \\ \text{or} \ q_w > X \ } } \log |\Delta_f F(P)|_w \right) \\ \hspace{-1cm} & \leq \left( \sum'_{\substack{w \in M_{K_i}, \\  q_w \leq X}} \sum_{P \in {\mathcal{S}}^{(i)} \setminus {\mathcal{S}}^{(i)}_F} \log |\Delta_f F(P)|_w \right)\,  + \, (\# {\mathcal{S}}^{(i)} \setminus {\mathcal{S}}^{(i)}_F \cdot \operatorname{max}_{P \in {\mathcal{S}}^{(i)}} h(\Delta_f F(P)).
\end{align*}
\medskip

\noindent
It follows from the construction of $F$ and earlier arguments that if $P \equiv O\  (\textrm{mod} \ w)$, we have $\log |\Delta_f F(P)|_w \leq -(T_0-T_f) \log q_w$, where $q_w$ is the norm of the prime ideal corresponding to $w$. Here $q_w$ is a power of prime lying outside $\mathcal{P}$. Reducing the N\'{e}ron model of $A/K$ at $w$, the reduced abelian variety $A/\F_{q_w}$ has cardinality $\sim q_w^g$ using the Hasse's bound. Now for a fix place $w$ and $a \in A(\F_{q_w})$, denote by $ {\mathcal{S}}^{(i)}_a$ the set consisting of elements of ${\mathcal{S}}^{(i)}$ that reduces to $a$ under reduction modulo $w$. This allows us to write

\begin{align}
\sum_{P \in {\mathcal{S}}^{(i)} \setminus {\mathcal{S}}^{(i)}_F} \sum_{w \in M_{K_i}} \log |F(P)|_w  \, \leq & \,-  (T_0-T_f) \sum'_{\substack{w \in M_{K_i}, \\  q_w \leq X}} \left( \sum_{a \in A(\mathbb{F}_{q_w})} (\#{\mathcal{S}}^{(i)}_a)^2- \#{\mathcal{S}}^{(i)} \cdot \#{\mathcal{S}}^{(i)}_F \right) \frac{\log q_w}{[K_i:\Q]} \nonumber \\ & + \, (\#{\mathcal{S}}^{(i)})^2 \cdot \operatorname{max}_{P \in {\mathcal{S}}^{(i)}} h(\Delta_f F(P)).
\end{align}
\medskip

\noindent
Since the LHS is zero by the product formula, we get 

\begin{align}
 (T_0-T_f) \sum'_{q \leq X} \left( \sum_{a \in A(\mathbb{F}_q)} (\#{\mathcal{S}}^{(i)}_a)^2- \#{\mathcal{S}}^{(i)} \cdot \#{\mathcal{S}}^{(i)}_F \right) \frac{\log q}{[K_i:\Q]} \nonumber \leq \, \#{\mathcal{S}}^{(i)} \cdot \operatorname{max}_{P \in {\mathcal{S}}^{(i)}} h(\Delta_f F(P)).
\end{align}
\medskip

\noindent
By the definition, $\sum_{a \in A(\mathbb{F}_q)} \, \#{\mathcal{S}}^{(i)}_a \, = \#{\mathcal{S}}^{(i)} $. Using the Cauchy–Schwarz inequality and $A(\mathbb{F}_q) \sim q^g$, we have 

\begin{align*}
(T_0-T_f) \sum'_{q \leq X } & \left( \frac{\#{\mathcal{S}}^{(i)}}{q^{2g}} - q^{2g} \ \#{\mathcal{S}}^{(i)} \cdot \# {\mathcal{S}}^{(i)}_F \right) \frac{N_q(K_i)}{[K_i : \Q]} \log q \hspace{1cm} \\ & \, \leq \,  (\# {\mathcal{S}}^{(i)})^2 \cdot \operatorname{max}_{P \in {\mathcal{S}}^{(i)}} h(\Delta_f F(P))\\
& \leq \, \, C_2 \cdot (\# {\mathcal{S}}^{(i)})^2 \cdot [h(F) + T_f \log(T_f + L) + T_f \log N + LN^2 \cdot \, \operatorname{max}_{P \in {\mathcal{S}}^{(i)}} \hat{h}(P)],
\end{align*}
where the last inequality follows from Theorem \ref{David} and Proposition \ref{Size of polynomial}.
\medskip

% \noindent
% Assuming $\#S(L)$ tends to infinity as $L$ is taken large, we obtain
% \begin{align*}
% T_0 \sum_q \psi_q \frac{\log q}{q^{2g}} \ll h(F) + L^2 \cdot \epsilon. 
% \end{align*}
% Now we combine all the places of norm $q=p^$. We denote by $N_q(L)$ the number of places of good reduction with norm $q$, where $q$ is a prime power. 
\noindent
Further, dividing both side by $T_0 \cdot (\# {\mathcal{S}}^{(i)})^2$ and using the fact that $\hat{h}(P) < 2 \epsilon$ for $P \in {\mathcal{S}}^{(i)}$, we get

\begin{align*}
\left(1-\frac{T_f}{T_0} \right) \sum'_{q \leq X} \left( \frac{1}{q^{2g}} - \frac{\# \mathcal{S}^{(i)}_F}{\# {\mathcal{S}}^{(i)}} \right) \frac{N_q(K_i)}{[K_i : \Q]} \log q &< \frac{C_2}{T_0} \cdot [h(F) + T_f \log(T_f + L) + T_f \log N \\ & + LN^2 \cdot \, 2\epsilon]
\end{align*}

\noindent
Substituting the values of parameters and taking $T_0$ large enough, we obtain

$$
 \sum'_{q \leq X} \left( \frac{1}{q^{2g}} - \frac{\# \mathcal{S}^{(i)}_F}{\# {\mathcal{S}}^{(i)}} \right) \frac{N_q(K_i)}{[K_i : \Q]} \log q < 2 C_2 \cdot \epsilon.
$$
\medskip

\noindent
Our assumption that $\{ P \in A(\mathcal{K})\} \ | \ \hat{h}(P) < \epsilon \}$ is an infinite set implies that $\# {\mathcal{S}}^{(i)} \rightarrow \infty $ as $i \rightarrow \infty$. This together with the fact that $\# {\mathcal{S}}_F^{(i)}$ remains bounded by $C_3 \cdot \left( \prod'_{q \leq X} q\right)^{2g^2} \cdot (\log T_0)^{4g}$ as we vary $i$, we deduce $\# \mathcal{S}^{(i)}_F/ \# {\mathcal{S}}^{(i)} \rightarrow 0$ as $i \rightarrow \infty$. Hence
$$
 \sum'_{q \leq X} \phi_q \frac{\log q}{q^{2g}} < 2\epsilon.
 $$

\noindent
 Thus, choosing $\epsilon = \frac{1}{2C_2}\sum'_{q \leq X} \phi_q \frac{\log q}{q^{2g}}$ gives a contradiction. Therefore, the set $\{ P \in A(\mathcal{K})\} \ | \ \hat{h}(P) < \epsilon \}$ is finite and hence
$$ 
\liminf_{P \in A(\mathcal{K})} \hat{h}(P) \geq \frac{1}{2C_2}\sum'_{q \leq X} \phi_q \frac{\log q}{q^{2g}}.
$$
\noindent
 Finally, we take $X \rightarrow \infty$ and replace $1/2C_2$ by $C$ to conclude the proof.
% \noindent
% Moreover, the additional constrain $\rho^g \log L \leq \frac{1}{2} \sum_{q \leq X} \psi_q \frac{\log q}{q^{2g}}$ gives

% $$
% \rho \sum'_{q \leq X} \frac{N_q(E)}{[E:\Q]} \frac{\log q}{q^{2g}} \leq \rho \frac{1}{2} \sum_{q \leq X} \psi_q \frac{\log q}{q^{2g}} + C_3 \epsilon
% $$
% \medskip

% \noindent
% or
% \begin{equation}\label{final}
% \rho \sum'_{q \leq X} \left(\frac{N_q(E)}{[E:\Q]} -\frac{1}{2} \psi_q  \right) \frac{\log q}{q^{2g}} \leq C_3 \epsilon.
% \end{equation}
% \medskip

% \noindent
% If the set $S'$ was infinite, we would've get a strictly increasing tower of number field $E_i$ such that $\# S(E_i) \rightarrow \infty$ and $E_i$ satisfying \eqref{final}. Since $\frac{N_q(E_i)}{[E_i:\Q]} \rightarrow \psi_q$ as $i \rightarrow \infty$, using \eqref{final}, we get

% $$
% \frac{\rho}{2} \sum'_{q \leq X} \psi_q  \frac{\log q}{q^{2g}} \leq C_3 \epsilon.
% $$

% \noindent
% Therefore, choosing $\epsilon \leq \frac{\rho}{C_3 2} \sum'_{q \leq X} \psi_q  \frac{\log q}{q^{2g}}$, we conclude that the set $S'$ is finite. Since $X$ was arbitrary, this conclude the proof of Theorem \ref{Main 2}.
\end{proof}

\section{\bf Proof of Theorem \ref{Main 3}}

The proof follows directly by combining Theorem \ref{Main 2} with a reduction argument of Baker and Silverman \cite[Proposition 2.1]{Baker-Silverman} for general abelian varieties and ample line bundles.

\begin{proof}[Proof of Theorem \ref{Main 3}]
 Let $A / K$ be an abelian variety of dimension $g$ and let $\rL$ be a symmetric ample line bundle on $A$. Over an algebraically closed field, Poincar\'{e}'s complete reducibility theorem tells that every abelian variety decomposes, up to isogeny, into a product of simple abelian varieties. We may thus find a finite extension $K^{\prime}$ of $K$ and geometrically simple abelian varieties $A_1, \ldots, A_r$ defined over $K^{\prime}$ so that there are isogenies
$$
\phi: A_1 \times \cdots \times A_r \longrightarrow A \quad \text { and } \quad \psi: A \longrightarrow A_1 \times \cdots \times A_r
$$
defined over $K^{\prime}$ with the property that $\phi \circ \psi=[m]$ for some integer $m \geq 1$.

\medskip
The line bundle $\phi^* \rL$ is ample on the product $A_1 \times \cdots \times A_r$ as $\phi$ is a finite. We fix an integer $n \geq 1$ so that $\phi^* \rL^{\otimes n}$ is very ample on the product, and we let $\rL_i^{\prime}=\left.\phi^* \rL^{\otimes n}\right|_{A_i}$ be the restriction to the $i^{\text {th }}$ factor. Then $\rL_i^{\prime}$ is a very ample line bundle on $A_i$.

\medskip
Let $P \in A(\bar{K})$ and write $\psi(P)=\left(P_1, \ldots, P_r\right)$. Then standard transformation properties of canonical height gives
$$
\begin{aligned}
\hat{h}_{A, \mathcal{L}}(P) & =\frac{1}{m^2} \hat{h}_{A, \mathcal{L}}([m] P) \\
& =\frac{1}{m^2 n} \hat{h}_{A, \mathcal{L}^{\otimes n}}([m] P) \\
& =\frac{1}{m^2 n} \hat{h}_{A_1 \times \cdots \times A_r, \phi^* \mathcal{L}^{\otimes n}}(\psi(P)) \\
& =\frac{1}{m^2 n}\left(\hat{h}_{A_1, \mathcal{L}_1^{\prime}}\left(P_1\right)+\hat{h}_{A_2, \mathcal{L}_2^{\prime}}\left(P_2\right)+\cdots+\hat{h}_{A_r, \mathcal{L}_r^{\prime}}\left(P_r\right)\right) .
\end{aligned}
$$

Let $\mathcal{K}$ be an infinite algebraic extension of $K$. Applying Theorem \ref{Main 2} to each abelian variety $A_i / K^{\prime}$ and very ample line bundle $\mathcal{L}_i^{\prime}$ to the infinite extension $\mathcal{K}K'$, we obtain
$$
\liminf_{P \in A_i(\mathcal{K}K')}\hat{h}_{A_i, \mathcal{L}_i^{\prime}}(P) \geq C(A, K', \mathcal{L}_i^{\prime}) \sum_{q=p^k}' \psi_q(\mathcal{K}K')\frac{\log q}{q^{2g}},
$$
where the sum runs over all prime powers $q=p^k$ with $p\notin \mathcal{P}$. Combining these estimates togethers, we get
$$
\liminf_{P \in A(\mathcal{K} K')}\hat{h}_{A, \mathcal{L}}(P) \geq \frac{1}{m^2 n} \left( \min _{1 \leq i \leq r} C(A, K', \mathcal{L}_i^{\prime}) \right) \sum_{q=p^k}' \psi_q(\mathcal{K}K')\frac{\log q}{q^{2g}},
$$
where the sum runs over all prime powers $q=p^k$ with $p\notin \mathcal{P}$. This completes the proof of Theorem \ref{Main 3}.

% \noindent
% Since $ P \in A_{\text {tors}} \Longleftrightarrow P_i \in\left(A_i\right)_{\text {tors }} \text { for all } 1 \leq i \leq r,$ we conclude that if $P \in A\left(\rK\right)$ is a non-torsion point, then
% $$
% \hat{h}_{A, \mathcal{L}}(P) \geq \frac{1}{m^2 n} \min _{1 \leq i \leq r} C_i .
% $$
% \noindent
\end{proof}

\section{Acknowledgement} I am grateful to Prof. Sinnou David for explaining the techniques used in \cite{David}, for many helpful discussions, and for providing crucial insights throughout this project. I am thankful to Prof. Aur\'{e}lien Galateau for explaining several technical details in \cite{Galateau} and for hosting me at the University of Cergy, where part of this work was carried out. I also thank Prof. Anup Dixit and Prof. Arnaud Plessis and Prof. Arno Fehm for valuable discussions and their comments on an earlier version of this manuscript.


\begin{thebibliography}{00}

\bibitem{Baker-2}
M. Baker, C. Petsche, Global discrepancy and small points on elliptic curves, \textit{Int. Math. Res. Not. (IMRN)}, \textbf{61}, pp. 3791-3834, (2005).

\bibitem{Baker-Silverman}
M. H. Baker, J.H. Silverman, A lower bound for the canonical height on abelian varieties over abelian extensions, \textit{Math. Res. Lett.}, \textbf{11}, pp. 377-396, (2004).

\bibitem{Bombieri-Gubler}
E. Bombieri, W. Gubler, Heights in Diophantine Geometry, New Math.
Monogr. \textbf{4}, \textit{Cambridge Univ. Press}, Cambridge, (2006).

\bibitem{BV}
E. Bombieri, J. Vaaler, On Siegel's lemma, \textit{Invent. Math.}, \textbf{73}, pp. 11-32, (1983).

\bibitem{BZ}
E. Bombieri, U. Zannier, A note on heights in certain infinite extensions of $\Q$, \textit{Atti. Acad. Naz. Lincei Cl. Sci. Mat. Fis. Natur. Rend. Lincei (9) Mat. Appl.}, \textbf{12}, pp. 5-14, (2001).

\bibitem{Carrizosa}
M. Carrizosa, Petits points et multiplication complexe,
\textit{Int. Math. Res. Not. (IMRN)}, \textbf{2009}, no. 16, pp. 3016–3097, (2009).

\bibitem{Checolli}
S. Checcoli, A. Fehm, On the Northcott property and local degrees, \textit{Proc. Amer. Math. Soc. }, \textbf{149}, pp. 2403–2414, (2021).

\bibitem{CL}
A. Chambert-Loir, Mesures et \'{E}quidistribution sur les espaces de Berkovich,
\textit{J. Reine Angew. Math.}, \textbf{595}, pp. 215–235, (2006).

\bibitem{David}
S. David, Fonctions th\^eta et points de torsion des vari\'et\'es ab\'eliennes, \textit{Compositio Math.}, {\bf 78}, no. 2, 121--160, (1991).

\bibitem{David-Hindry}
S. David, M. Hindry, Minoration de la hauteur de N\'{e}ron-Tate sur les vari\'{e}t\'{e}s ab\'{e}liennes de type C. M., \textit{J. Reine Angew. Math.}, no. \textbf{529}, pp. 1-74, (2000).

\bibitem{DP}
S. David, P. Philippon, Minorations des hauteurs normalis\'ees des sous-vari\'et\'es des tores, \textit{Ann. Scuola Norm. Sup. Pisa Cl. Sci. (4)}, {\bf 28}, no. 3, pp. 489-543, (1999).

\bibitem{Dimitrov}
V. Dimitrov, \hyperlink{https://mathoverflow.net/questions/123096}{https://mathoverflow.net/questions/123096}, (2013).

\bibitem{AB-SK}
A. B. Dixit, S. Kala, Bogomolov property for certain infinite non-Galois extensions, \textit{Ann. Scuola Norm. Sup. Pisa Cl. Sci.}, {\bf 31}, (2026).

\bibitem{Fili1}
P. Fili, On the heights of totally p-adic numbers, \textit{J. Th\'{e}or. Nombres Bordeaux}, \textbf{26}(1), pp. 103-109, (2014).

\bibitem{Galateau}
A. Galateau, Le probl\`eme de Bogomolov effectif sur les vari\'et\'es ab\'eliennes, \textit{Algebra Number Theory}, {\bf 4}, no.~5, pp. 547-598, (2010).

\bibitem{Gubler}
W. Gubler, The Bogomolov conjecture for totally degenerate abelian varieties, \textit{Invent. math.}, \textbf{169}, pp. 377–400, (2007). 

\bibitem{Hindry-Silverman}
M. Hindry , J. H. Silverman, Diophantine Geometry : An Introduction, Graduate Texts in Mathematics \textbf{201}, \textit{Springer-Verlag}, New York, pp. xiv+558, (2000).

\bibitem{AbelMasser}
 D. Masser, Small values of the quadratic part of the N\'{e}ron-Tate height on an abelian variety, \textit{Compositio Mathematica}, \textbf{53}, pp. 153-170, (1984).

 \bibitem{Philippon}
 P. Philippon, Lemmes de z\'{e}ros dans les groupes alg\'{e}briques commutatifs. Bulletin de la Soci\'{e}t\'{e} Math\'{e}matique de France, Vol. \textbf{114}, pp. 355-383, (1986).

 \bibitem{Plessis}
 A. Plessis and S. Sahoo, Lower bounds for heights on some algebraic dynamical systems, \textit{arXiv:2502.03039}, (2025).

 \bibitem{Roy}
 D. Roy and J. L. Thunder, An absolute Siegel's lemma, \textit{J. Reine Angew. Math.}, {\bf 476}, pp. 1-26, (1996).

\bibitem{Schmidt} W. M. Schmidt, { Diophantine approximations and Diophantine equations}, Lecture Notes in Mathematics, 1467, \textit{Springer}, Berlin, (1991).

 % \bibitem{AWS}
 % J. H. Silverman, Canonical Heights on Abelian Varieties, \textit{Lecture Notes for the Arizona Winter School}, March 2–6, 2024, \hyperlink{https://swc-math.github.io/aws/2024/2024SilvermanNotes.pdf}{https://swc-math.github.io/aws/2024/2024SilvermanNotes.pdf}, (2024). 

 \bibitem{SUZ}
 L. Szpiro, E. Ullmo, and S. W. Zhang, \'{E}quir\'{e}partition des petits points, \textit{Invent. Math.}, \textbf{127}, no. 2, pp. 337–347, (1997).

% \bibitem{Zarhin}
%  Yu. G. Zarhin, Yu. I. Manin, Height on families of Abelian varieties, \textit{Math. USSR-Sb.}, \textbf{18}, 2, pp. 169–179, (1972).

 \bibitem{Zhang} S. W. Zhang, Equidistribution of small points on abelian varieties, \textit{Ann. of Math. (2)}, \textbf{147}, no. 1, pp. 159-165, (1998).

\end{thebibliography}
\end{document}